\def \version {October 2, 2025}

\documentclass[12pt]{article}

\usepackage{amsmath}
\usepackage{amssymb}
\usepackage{amsthm}
\usepackage{stmaryrd}
\usepackage{graphicx} 
\usepackage{subcaption}
\usepackage{tikz}
\usepackage{lineno}
\usepackage{comment}

\newcommand{\komm}[1]{{\bsk \vbox{\hrule \ssk {\nin \large\sf #1} \ssk \hrule \par} \bsk }}

\newcommand{\nev}[1]{{\bf\itshape (#1)}\, }

\def \ssk {\smallskip}
\def \msk {\medskip}
\def \bsk {\bigskip}
\def \nin {\noindent}
\def \qed {\hfill $\Box$}

\newcommand{\floor}[1]{\lfloor #1 \rfloor }
\newcommand{\ceil}[1]{\lceil #1 \rceil }
\newtheorem{Theorem}{Theorem}
\def \btm {\begin{Theorem}}
\def \etm {\end{Theorem}}
\newtheorem{Lemma}[Theorem]{Lemma}
\def \blm {\begin{Lemma}}
\def \elm {\end{Lemma}}
\newtheorem{Problem}[Theorem]{Problem}
\def \bpm {\begin{Problem}}
\def \epm {\end{Problem}}
\newtheorem{Proposition}[Theorem]{Proposition}
\def \bpn {\begin{Proposition}}
\def \epn {\end{Proposition}}
\newtheorem{Corollary}[Theorem]{Corollary}
\def \bcr {\begin{Corollary}}
\def \ecr {\end{Corollary}}
\newtheorem{con}[Theorem]{Conjecture}
\def \bcj {\begin{con}}
\def \ecj {\end{con}}
\newtheorem{Homework}[Theorem]{Homework}
\def \bhw {\begin{Homework}}
\def \ehw {\end{Homework}}
\newtheorem{Example}[Theorem]{Example}
\def \bex {\begin{Example}\rm }
\def \eex {\end{Example}}
\newtheorem{Definition}[Theorem]{Definition}
\def \bdf {\begin{Definition}\rm }
\def \edf {\end{Definition}}
\newtheorem{Remark}[Theorem]{Remark}
\def \brm {\begin{Remark}\rm }
\def \erm {\end{Remark}}
\def \bpf {\begin{proof}}
\def \epf {\end{proof}}
\def \hsp {\hspace{0.5em}}

\def \Boxx {\,\Box\,}
\def \strg {\,\boxtimes\,}
\def \nnn {\mathbb{N}}
\def \zzz {\mathbb{Z}}

\def \soc {\chi_\mathrm{so}}
\def \sid {\alpha_\mathrm{od}}
\def \siid {\alpha_\mathrm{iod}}
\def \rhodd {\varrho_\mathrm{od}}

\def \vp {\varphi}
\def \es {\varnothing}
\def \smin {\,\diagdown\,}
\def \soi {odd independent}
\def \sois {\soi\ set}
\def \soin  {odd independence number}

\def \sokk  {strong odd coloring}
\def \gdinf {G_{d\,;\,\infty}}
\def \KG {\mbox{\sf KG}}

\newcommand{\Mod}[1]{\ [\mathrm{mod}\ #1]}

\begin{document}


\title{The odd independence number of graphs, II: Finite and infinite grids and chessboard graphs}
\author{Yair Caro\thanks{\hsp Department of Mathematics,
  University of Haifa-Oranim, Tivon 36006, Israel}\,,
 Mirko Petru\v sevski\thanks{\hsp Faculty of Mechanical Engineering, Ss. Cyril and Methodius University in Skopje, Macedonia}\,,
  Riste \v Skrekovski\thanks{\hsp Faculty of Mathematics and Physics, University of Ljubljana;  Faculty of Information Studies in Novo Mesto;  Rudolfovo - Science and Technology Centre Novo Mesto, FAMNIT, University of Primorska, Slovenia}\,,
   Zsolt Tuza\thanks{\hsp Alfr\'ed R\'enyi Institute of Mathematics,
    H-1053 Budapest, Re\'altanoda u.~13--15, Hungary;
     and
    Department of Computer Science and Systems Technology,
    University of Pannonia, 8200 Veszpr\'em, Egyetem u.~10, Hungary}}
\date{\small Latest update on \version}
\maketitle

\begin{abstract}

An odd independent set $S$ in a graph $G=(V,E)$ is an independent set of
 vertices such that, for every vertex $v \in V \smin S$, either
 $N(v) \cap S = \es$  or  $|N(v) \cap S| \equiv 1$ (mod 2),
 where $N(v)$ stands for the open neighborhood of $v$.
The largest cardinality of odd independent sets of a graph $G$,
  denoted  $\sid(G)$, is called the \soin\ of $G$.

This new parameter is a natural companion to the recently introduced
 strong odd chromatic number.
A proper vertex coloring of a graph $G$ is a strong odd coloring
 if, for every vertex $v \in V(G)$, each color used in the neighborhood
 of $v$ appears an odd number of times in $N(v)$.
The minimum number of colors in a strong odd coloring of $G$
 is denoted by $\soc(G)$.

A simple relation involving these two parameters and the order $|G|$
 of $G$ is $\sid(G)\cdot\soc(G) \geq |G|$, parallel to the same on
 chromatic number and independence number.  

In the present work, which is a companion to our first paper on the subject [The odd independence number of graphs, I: Foundations and classical classes],
we focus on grid-like and chessboard-like graphs and compute or estimate their odd independence number and their strong odd chromatic number.
Among the many results obtained, the following give the flavour of this paper:

(1)\quad $0.375 \leq \rhodd(P_\infty \Boxx P_\infty) \leq 0.384615...$, where $\rhodd(P_\infty \Boxx P_\infty)$ is the odd independence ratio.

(2)\quad $\soc(\gdinf) = 3$ for all $d \geq 1$, where $\gdinf$ is the infinite $d$-dimensional grid. As a consequence, $\rhodd(\gdinf) \geq 1/3$.

(3)\quad The $r$-King graph $G$ on $n^2$ vertices has $\sid(G) = \ceil{n/(2r+1)}^2$. Moreover, $\soc(G) = (2r + 1)^2$ if $n \geq 2r + 1$, and $\soc(G) = n^2$ if $n \leq 2r$.

Many open problems are given for future research.

\bsk

\nin

\bsk

\nin

\end{abstract}

\msk

\section{Introduction}

An \emph{odd independent set} $S$ in a graph $G = (V, E)$ is an independent set of vertices such that, for every vertex $v \in V \smin S$, either $N (v) \cap S = \es$ or $|N (v) \cap S| \equiv 1$ (mod 2), where $N (v)$ stands for the open neighbourhood  of $v$.
The largest cardinality of odd independent sets in $G$, denoted $\sid(G)$, is called the \emph{odd independence number} of $G$.
This new parameter is a natural companion to the recently introduced strong odd chromatic number \cite{CPST,GK-etal,KP-a24,PMF-DM26,P-a25}.
A proper vertex coloring of a graph $G$ is a strong odd coloring if, for every vertex $v \in V (G)$, each
color used in the neighborhood of $v$ appears an odd number of times in $N (v)$. The minimum number of colors in a strong odd coloring of $G$ is denoted by $\soc(G)$.

After the initial paper \cite{part-1} introducing and starting a systematic study of the odd independence number, the present work is aimed to consider in depth the odd independence number of various types of grids and chessboard-like graphs, both finite and infinite.

In Section \ref{s:basic} we present basic upper bounds for the odd independence number of regular or nearly regular $K_{1,r}$-free graphs.
These will be then applied for the various grids which are $K_{1,3}$-free or $K_{1,5}$-free.

In section \ref{s:grid} we consider the odd independence number of three types of grids: $\sid(P_k \Boxx P_k)$, $\sid(P_k \Boxx C_k)$, and $\sid(C_k \Boxx C_k)$. Further, we also study the infinite grid $P_\infty \Boxx P_\infty$. After some technical preparation we prove
that the maximum density $\rhodd(P_\infty \Boxx P_\infty)$ of odd independent sets in this grid can be bounded as $3/8 = 0.375 \leq \rhodd(P_\infty \Boxx P_\infty) \leq 5/13 < 0.38462$.

For the infinite $d$-dimensional grid $\gdinf := P_\infty \Boxx \cdots \Boxx P_\infty$ on the vertex set $\zzz^d$ we have $\chi((\gdinf)^2) = 2d + 1$.
In sharp contrast to this, the strong odd chromatic number of the $d$-dimensional infinite grid is independent of the dimension, and we prove
$\soc(\gdinf) = 3$ for all $d \geq 1$. As a consequence, the maximum density of odd independent sets in these grids satisfies $\rhodd(\gdinf) \geq 1/3$.

In section \ref{s:chess} we consider various chessboard-like grid graphs.
Among the many results proved, we mention here the $r$-King graph.
The $r$-King graph, $\KG_{n,r}$ ($r \in \nnn$) of order $n \times n$, has the same vertex set as $P_n \Boxx P_n$.
Two vertices $(a,b),(a',b')\in\{1,\dots,n\}\times\{1,\dots,n\}$ are adjacent if and only if $|a-a'|\leq r$ and $|b-b'|\leq r$.
The case $r = 1$ is the classical King graph.
We prove that the $r$-King graph of order $n \times n$ has $\sid(\KG_{n,r}) = \ceil{n/(2r+1)}^2$. Moreover, $\soc(\KG_{n,r}) = (2r + 1)^2$ if $n \geq 2r + 1$, and $\soc(\KG_{n,r}) = n^2$ if $n \leq 2r$.
Similarly, we develop lower and upper bounds for the odd independence numbers of $r$-Rook graphs, $r$-Bishop graphs, and their infinite graph version, as well as the Knight graph and the Queen graph.

In section \ref{s:tri-hex} we consider the infinite triangular grid $T_\infty$, proving that the maximum density of odd independent sets is $\rhodd (T_\infty)= 1/3$; also, $\soc(T_\infty)= 3$. Moreover, for the infinite hexagonal grid $H_\infty$ we prove
$\rhodd(H_\infty) = 1/2$ and $\soc(H_\infty) = 2$.

In section \ref{s:concl} we offer many open problems and conjectures concerning odd independence in various grid and chessboard-like graphs, as well as concerning density of odd independence in the infinite grid and chessboard graphs.

The Appendix contains data obtained with computer and used in the various tables in this paper.

\subsection{Notation}

For a graph $G$, standard notation inclues
 the order $|G|:=|V(G)|$ (the number of vertices in $G$),
 the independence number $\alpha(G)$,
 the chromatic number~$\chi(G)$,
 and the square $G^2$, in which
 two vertices are adjacent if and
   only if they are at distance at most 2 apart in $G$ (i.e., they are
   adjacent or have a common neighbor in $G$).

The open and closed neighborhood of vertex $v$ is denoted
 by $N(v)$ and $N[v]$, respectively.
A graph is said to be $K_{1,r}$-free ($r\geq 3$) if it contains
 no \emph{induced} subgraph isomorphic to $K_{1,r}$.
The common term for $K_{1,3}$-free is claw-free.

As a less common notation, we write $p$ [mod $q$] for
 the residue of $p$ modulo $q$ (where $p\geq 0$ and $q\geq 2$ are integers).

\section{Some basic tools}
\label{s:basic}

This section in mainly about obtaining upper bounds on the odd independence number for $d$-regular or nearly $d$-regular $K_{1,r}$-free graphs.
These estimates are then used to upper-bound $\sid$ in various grid and chessboard-like graphs which are nearly regular and either $K_{1,3}$-free or $K_{1,5}$-free.

\bpn
\label{p:reg-r-claw}
Let $G$ be a $d$-regular, $K_{1,r}$-free graph ($r \geq 3$) on $n$ vertices.

 \begin{itemize}
  \item[$(i)$] If $r \equiv 0$ {\rm (mod 2)}, then
   $\sid(G) \leq \frac{r-1}{d+r-1}\,n$\,.
  \item[$(ii)$] If $r \equiv 1$ {\rm (mod 2)}, then
   $\sid(G) \leq \frac{r-2}{d+r-2}\,n$\,.
 \end{itemize}
\epn

\bpf
Let $S$ be a maximum odd independent set in $G$.
Consider any $v \in V \smin S$.
Either $v$ is not adjacent to any vertex in $S$, or it has at least one neighbor $u \in S$.
In the latter case it can have at most $r - 1$ neighbors in $S$ if $r$ is even, and at most $r - 2$ neighbors if $r$ is odd,
 since $r$ independent neighbors would yield an induced $K_{1,r}$ subgraph.
 Now count the number of edges between $S$ and $V \smin S$ in two ways and obtain
 $d\cdot |S| \leq (r-1)(n-|S|)$ if $r$ is even, and
 $d\cdot |S| \leq (r-2)(n-|S|)$ if $r$ is odd.
Hence it follows that $|S| \leq (r - 1)n/(d + r -1)$ or $|S| \leq (r - 2)n/(d + r -2)$, respectively.
\epf

The above estimates are tight for infinitely many
 combinations of the parameters, as shown next.

\bex
$(i)$\quad
Let $p\geq 3$ be odd, and $t\geq 2$ any integer.
Consider the complete $t$-partite graph $G=K_{t*p}$,
 with $t$ vertex classes of size $p$ each.
This $G$ is $K_{1,r}$-free with $r=p+1$ even.
The degree of regularity is $d=(t-1)\cdot p$, the order is $n=pt$, and the
 obtained upper bound is $(r-1)\cdot n/(d +r - 1) = tp^2/((t-1)\cdot p + p) = p$\,;
 and each vertex class is indeed an odd independent set of cardinality $p$.
This shows the tightness of $(i)$ for all $r\geq 4$ even, and all $d\equiv 0$ {\rm (mod $r-1$)}.

$(ii)$\quad
Taking $G=K_{t*p}$ with $p\geq 2$ even, the graph is
 $K_{1,r}$-free with $r=p+1$ odd.
In this case we have the upper bound
 $(r-2)\cdot n/(d +r - 2) = tp\cdot(p-1)/((t-1)\cdot p + p-1) = tp\cdot(p-1)/(tp-1) = p-1 + (p-1)/(tp-1)$,
 whose integer part is $p-1$ for all $t\geq 2$.
This is tight for $(ii)$, because the largest odd independent sets are
 obtained by omitting one vertex from any one vertex class, hence
 of cardinality $p-1$.
\eex

Also, the following asymptotic version of the above proposition is useful.

\bpn
\label{p:reg-clawfree-asymp}
Let $\Delta$ be a natural number, and for $j=1,2,\dots$ let
 $G_j$ be a $K_{1,r}$-free graph ($r\geq 3$) on $n_j$ vertices, $n_j \to\infty$,
 of maximum degree $\Delta$ in which the number $z_j$ of vertices
 of degrees smaller than $\Delta$ satisfies $z_j = o(n_j)$.

 \begin{itemize}
  \item[$(i)$] If $r \equiv 0$ {\rm (mod 2)}, then
   $\limsup \sid(G_j) / n \leq \frac{r-1}{\Delta+r-1}$\,.
  \item[$(ii)$] If $r \equiv 1$ {\rm (mod 2)}, then
   $\limsup \sid(G_j) / n \leq \frac{r-2}{\Delta+r-2}$\,.
 \end{itemize}
\epn

\bpf
As in the the proof of Proposition \ref{p:reg-r-claw},
 for the size of a maximum odd-independent set $S$ we obtain the inequality
 $d\cdot (|S|-z_j) \leq (r-1)(n_j-|S|)$ or $d\cdot (|S|-z_j) \leq (r-2)(n_j-|S|)$,
 depending on whether $r$ is even or odd.
Since $z_j=o(n_j)$, the asymptotic upper bounds follow.
\epf

Among the chessboard graphs studied in Section \ref{s:chess} below,
 Bishop and Rook are $K_{1,3}$-free,
 whereas
  Queen, $r$-Rook and $r$-Bishop (and also the grid graphs) are $K_{1,5}$-free,
 all with even maximum degree.
For $K_{1,3}$-free boards we know directly $\sid(G) = \alpha(G^2)$, and for
 $K_{1,5}$-free boards we obtain the upper-bound ratio of $3/(\Delta+3)$.

The following concept of forbidden pair / forcing pair, introduced in Part~1 \cite{part-1}, will be also useful here.

\bdf   \nev{Forbidden pair \& Forcing pair}
\

\begin{itemize}
\item 
Two nonadjacent vertices $x,y$ form a \emph{forbidden pair} if they have a common neighbor $z$ such that $N[z]\subseteq N[x]\cup N[y]$.

\item 
Two nonadjacent vertices $x,y$ form a \emph{forcing pair} if they have a common neighbor $z$ such that, for each vertex $w$ in $N(v)$ which is independent of $x$ or $y$ or both, either $x,w$ or $y,w$ or both are forbidden pairs.
\end{itemize}
\edf

\blm[\cite{part-1}]
\label{l:force}
If $S\subseteq V(G)$ is an odd-independent set, then it contains no forbidden pair, nor a forcing pair.
\elm

\section{Three types of square grids}
\label{s:grid}

Here we consider grids $P_k \Boxx P_k$ in the plane, $P_k \Boxx C_k$
 on the cylinder, and $C_k \Boxx C_k$ on the torus.
Moreover, the infinite grid in the plane will be investigated, and a
 strong odd coloring of higher-dimensional grids is also provided.

Before all these, let us state that the odd independence number is
 supermultiplicative with respect to the Cartesian product operation.

\bpn
For any two graphs $G$ and $H$ the general lower bound
 $\sid(G \Boxx H)\geq \sid(G) \cdot \sid(H)$ is valid.
\epn

\bpf
Let $S_G = \{ u_1,\dots,u_g \}$ and $S_H = \{ v_1,\dots,v_h \}$
 be maximum odd independent sets of $G$ and $H$, respectively
 (where $g=\sid(G)$ and $h=\sid(H)$).
We claim that the set $S_G \times S_H = \{ (u_i,v_j) \mid
 1\leq i\leq g, \ 1\leq j\leq h \}$ is odd-independent in $G \Boxx H$.

The graph $G \Boxx H$ is the edge-disjoint union of
 $|H|$ copies of $G$ and $|G|$ copies of $H$.
More specifially, each vertex $w=(u_a,v_b) \in G \Boxx H$
 is in one copy $G_b$ of $G$ and in one copy $H_a$ of $H$.
This immediately implies that $S_G \times S_H$ is independent.
Assuming that $w=(u_a,v_b) \notin S_G \Boxx S_H$,
 we claim that if $w$ has at least one neighbor
 (and then by assumption an odd number of them)
 in $S_G\cap G_b$, then it has no neighbor in $S_H\cap H_a$,
 and vice versa.
Inded, suppose that both $(u_i,v_b),(u_a,v_j)\in
 S_G \times S_H$ are adjacent to
 $(u_a,v_b)$.
This can occur only if both $u_a,u_i\in S_G$.
But then adjacency of $(u_a,v_b)$ and $(u_i,v_b)$ implies
 $u_iu_a\in E(G)$, contradicting the assumption that
 $S_G$ is independent.
\epf

Products of complete graphs establish infinitely many cases of equality.
However, $\sid(C_3 \Boxx C_4) = 4$ while
 $\sid(C_3) = \sid(C_4) = 1$.
Hence already
 very small examples show that multiplicativity
 with respect to Cartesian product does not hold in general, not even in cases where
 one of the factors is a complete graph.

\begin{table}[ht]
\begin{center}
    \begin{tabular}{c|cccccccccccc}
       $k$  & 3 & 4 & 5 & 6 & 7 & 8 & 9 & 10 & 11 & 12 & 13 & 14   \\
\hline
      $P_k \Boxx P_k$   & 5 & 5 & 12 & 12 & 20 & 21 & 29 & 33 & 42 & 48 & 60 & 64  \\
      $P_k \Boxx C_k$   & 2 & 5 & 6 & 12 & 14 & 24 & 25 & 34 & 37 & 48 & 52 & 70  \\
      $C_k \Boxx C_k$   & 1 & 6 & 5 & 12 & 12 & 24 & 21 & 30 & 34 & 54 & 49 & 62 
    \end{tabular}
    \caption{Values of $\sid$ in small square grids}
    \label{tab:grid}
\end{center}
\end{table}

\bpn
For $k\leq 14$, the values $\sid(P_k \Boxx P_k)$, $\sid(P_k \Boxx C_k)$, and $\sid(C_k \Boxx C_k)$ are as shown in Table \ref{tab:grid}.
\epn

\bpf
These data have been determined using computer search.
The lower bounds are obtained by explicit constructions, which are listed in the appendix.
\epf

Our next concern is the proportion of vertices in an \sois.

\bdf
If $G$ is a finite graph, we denote $\rhodd(G) := \sid(G)/|G|$ and
 call it the \emph{odd independence ratio} of $G$.
Further, if $G$ is an infinite plane graph with a given planar embedding,
 we define the \emph{odd independence density} of $G$ as $\rhodd(G) :=
  {\displaystyle \sup_S \, 	\limsup_{H_n\,, \ n\to\infty} } |S\cap V(H_n)|/|H_n|$
  taken over all \sois s $S$ of $G$ and all $H_n$ induced by the
 vertices inside axis-parallel squares of size $n$ in the plane.
\edf

\bdf
Given $k\geq 3$, an \emph{internally odd-independent set} is an independent set
 $S\subset V(P_k\Boxx P_k)$ such that every vertex of degree 4 in $P_k\Boxx P_k$
 (i.e., those not on the boundary) either belongs to $S$ or has an
 odd number of neighbors in $S$.
The largest possible $|S|$ of this kind is denoted by $\siid(P_k\Boxx P_k)$.
\edf

For $k<3$ there are no internal vertices in $P_k\Boxx P_k$, this is the reason why we start the definition from 3.
As we shall see in Theorem \ref{t:rho-grid}, the notion of internally odd independence is useful in upper-bounding the density of odd independence in the infinite 2-dimensional grid.

\bpn
The values of $\siid(P_k\Boxx P_k)$ for $k\leq 26$ are equal to those given in table \ref{tab:intern-odd}. 
\epn

\bpf
These values have been determined via exhaustive computer search.
\epf

\begin{table}[ht]
\begin{center}
{
\scriptsize
    \begin{tabular}{ccccccccc}
       $k$  & 3 & 4 & 5 & 6 & 7 & 8 & 9 & 10  \\
\hline
      $\siid$   & 5 & 7 & 12 & 15 & 21 & 26 & 34 & 40  \\
      $\siid/k^2$   & 0.555556 & \textbf{0.4375} & 0.48 & \textbf{0.416667} & 0.428571 & \textbf{0.40625} & 0.419753 & \textbf{0.4}  \\
      &&&&&&&& \\
       $k$  & 11 & 12 & 13 & 14 & 15 & 16 & 17 & 18  \\
\hline
      $\siid$   & 49 & 57 & 68 & 77 & 89 & 100 & 114 & 125  \\
      $\siid/k^2$   & 0.404959 & \textbf{0.395833} & 0.402367 & \textbf{0.392857} & 0.395556 & \textbf{0.390625} & 0.394464 & \textbf{0.388889}  \\
      &&&&&&&& \\
       $k$  & 19 & 20 & 21 & 22 & 23 & 24 & 25 & 26  \\
\hline
      $\siid$   & 141 & 155 & 172 & 187 & 204 & 222 & 242 & 260  \\
      $\siid/k^2$   & 0.390582 & \textbf{0.3875} & 0.390023 & \textbf{0.386364} & 0.387524 & \textbf{0.385417} & 0.3872 & \textbf{0.384615}  
    \end{tabular}
    \caption{Values of $\siid(P_k\Boxx P_k)$, exact or rounded to six decimal digits}
    \label{tab:intern-odd}
}
\end{center}
\end{table}

One can observe that the subsequences of fixed parity are monotone
 decreasing, for both $k$ odd and $k$ even.
However, the value of $\siid/k^2$ for $k=2t-1$ and $k=2t+1$ are greater
 than that for $k=2t$, in all of the investigated cases.
It is our belief that this property remains valid for all $k$, and that
 the limit of the sequence is equal to $3/8$.

\begin{figure}[ht]
\begin{center}
\begin{tikzpicture}
\draw[step=1.5cm,black,ultra thick] (0,0) grid (4.5,4.5);
\fill[orange] (0,4.5) circle(6pt);
\fill[orange] (3,4.5) circle(6pt);
\fill[orange] (1.5,3) circle(6pt);
\fill[orange] (0,1.5) circle(6pt);
\fill[orange] (3,1.5) circle(6pt);
\fill[orange] (4.5,0) circle(6pt);
\end{tikzpicture}
\caption{$\sid(C_4\Boxx C_4)=6$.}
\label{fig:4x4}
\end{center}
\end{figure}

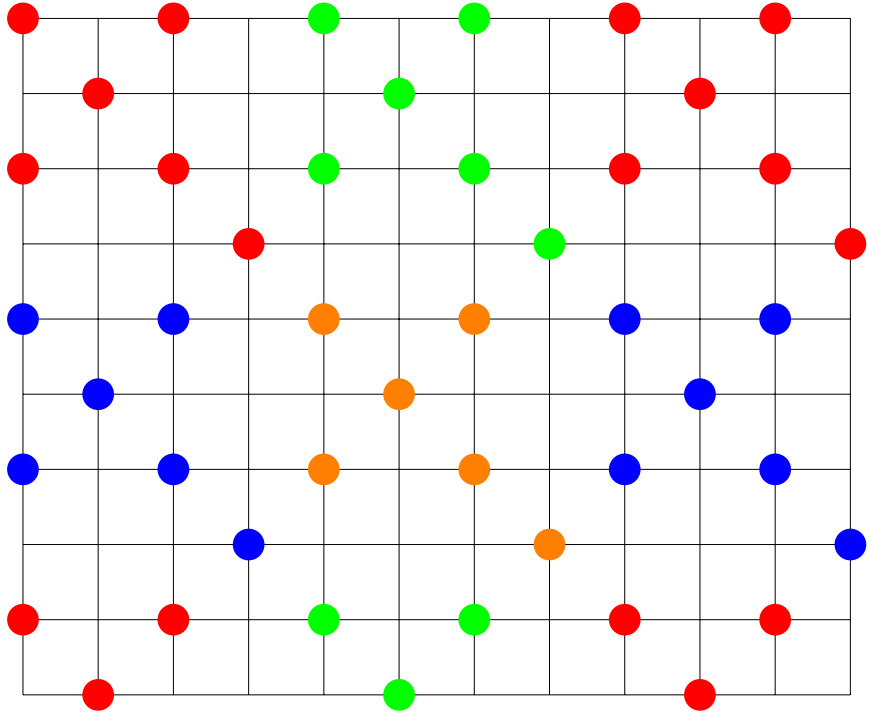
\begin{figure}[ht]
\begin{center}
\begin{tikzpicture}
\draw[step=1cm,black,very thin] (0,0) grid (11,9);
\fill[orange] (4,5) circle(6pt);
\fill[orange] (6,5) circle(6pt);
\fill[orange] (5,4) circle(6pt);
\fill[orange] (4,3) circle(6pt);
\fill[orange] (6,3) circle(6pt);
\fill[orange] (7,2) circle(6pt);

\fill[blue] (8,5) circle(6pt);
\fill[blue] (10,5) circle(6pt);
\fill[blue] (9,4) circle(6pt);
\fill[blue] (8,3) circle(6pt);
\fill[blue] (10,3) circle(6pt);
\fill[blue] (11,2) circle(6pt);

\fill[blue] (0,5) circle(6pt);
\fill[blue] (2,5) circle(6pt);
\fill[blue] (1,4) circle(6pt);
\fill[blue] (0,3) circle(6pt);
\fill[blue] (2,3) circle(6pt);
\fill[blue] (3,2) circle(6pt);

\fill[red] (0,1) circle(6pt);
\fill[red] (2,1) circle(6pt);
\fill[red] (1,0) circle(6pt);
\fill[red] (0,7) circle(6pt);
\fill[red] (2,7) circle(6pt);
\fill[red] (3,6) circle(6pt);

\fill[red] (8,1) circle(6pt);
\fill[red] (10,1) circle(6pt);
\fill[red] (9,0) circle(6pt);
\fill[red] (8,7) circle(6pt);
\fill[red] (10,7) circle(6pt);
\fill[red] (11,6) circle(6pt);

\fill[green] (4,1) circle(6pt);
\fill[green] (6,1) circle(6pt);
\fill[green] (5,0) circle(6pt);
\fill[green] (4,7) circle(6pt);
\fill[green] (6,7) circle(6pt);
\fill[green] (7,6) circle(6pt);

\fill[red] (1,8) circle(6pt);
\fill[green] (5,8) circle(6pt);
\fill[red] (9,8) circle(6pt);

\fill[red] (0,9) circle(6pt);
\fill[red] (2,9) circle(6pt);
\fill[green] (4,9) circle(6pt);
\fill[green] (6,9) circle(6pt);
\fill[red] (8,9) circle(6pt);
\fill[red] (10,9) circle(6pt);

\end{tikzpicture}
\caption{Detail of independent set of density $3/8$ in $P_\infty\Boxx P_\infty$, generated from $\sid(C_4\Boxx C_4)=6$.}
\label{fig:plane-LB-3/8}
\end{center}
\end{figure}

\btm
\label{t:rho-grid}
${\displaystyle \sup_{k\geq 3} } \ \rhodd(C_k\Boxx C_k) \leq \rhodd(P_\infty\Boxx P_\infty) \leq \inf \siid(P_k\Boxx P_k)/k^2$.
\etm

\begin{figure}[ht]
\begin{center}
\resizebox{5cm}{!}{%
\includegraphics
{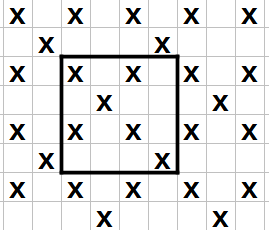}
  }
\caption{Cellular representation, also indicating the generating pattern which yields $\rhodd\geq 3/8$.}
\label{fig:grid-cells}
\end{center}
\end{figure}

\bpf
For the lower bound we partition the plane into $k\times k$ squares.
More explicitly, we copy an optimal \sois\ $S_k$ of $C_k\Boxx C_k$
 inside the integer points $\{1,\dots, k\} \times \{1,\dots, k\}$
 and take its translates by the vectors $(ak,bk)$ for all $a,b\in\zzz$;
 see Figure \ref{fig:4x4} with a pattern for $k=4$, and a detail of the
 corresponding \sois\ of the infinite grid exhibited in
  Figure~\ref{fig:plane-LB-3/8}.
In this way starting from $C_k\Boxx C_k$ we obtain an \sois\ of
 density $\sid(C_k\Boxx C_k)/k^2$, hence a lower bound.
The additional Figure~\ref{fig:grid-cells} represents the vertices of the grid by
 small square cells.
The cells corresponding to the vertices of
 the selected \sois\ $S$ are marked with {\Large\sf x}.
Further, thick frame surrounds the pattern whose translates yield $S$.
This $4\times 4$ pattern is the equivalent of the set $S_4$
 exhibited in Figure \ref{fig:4x4}.

For the upper bound, partition the infinite grid into vertex-disjoint
 copies of $P_k \Boxx P_k$.
If $S$ is an \sois\ of the grid, by definition it can have at most
 $\siid(P_k \Boxx P_k)$ vertices in each copy.
Consider a sequence $H_n \cong P_n \Boxx P_n$, $n\to\infty$,
 establishing $\rhodd(G) := \sup_S {\displaystyle\lim_{
   n\to\infty} } |S\cap V(H_n)|/|H_n|$.
Writing $n$ in the form $n=qk+r$ ($0\leq r\leq k-1$), we pack $q^2$
 vertex-disjoint copies of $P_k \Boxx P_k$ in $H_n$ and obtain

 $$
   \rhodd(G) = \sup_S {\lim_{ n\to\infty} } |S\cap V(H_n)|/|H_n| \leq
    \frac{q^2\cdot \siid(P_k \Boxx P_k) + n^2-(qk)^2}{n^2}
 $$

 $$
   \leq \frac{\siid(P_k \Boxx P_k)}{k^2} + \frac{2rn-r^2}{n^2}
    < \frac{\siid(P_k \Boxx P_k)}{k^2} + \frac{2k}{n} \,.
 $$
The last term tends to 0 for any $k$ as $n\to\infty$.
Thus, the density of $S$ is at most $\siid(P_k \Boxx P_k)/k^2$.
Taking a sequence of $k$ that attains the infimum of this ratio,
 the theorem follows.
\epf

\bcr
$\sid(C_4\Boxx C_4)/4^2 \leq \rhodd(P_\infty\Boxx P_\infty) \leq \siid(P_{26}\Boxx P_{26})/26^2$; more explicitly, $\rhodd(P_\infty\Boxx P_\infty)$ is between $3/8 =0.375$ and $260/676 = 5/13 \approx 0.38461538461538464$.
\ecr

This result substantially improves the upper bound of 3/7
 obtained from Proposition \ref{p:reg-clawfree-asymp} using the fact that
 the grid graphs are $K_{1,5}$-free and have maximum degree 4.

We close this subsection with results on the
 infinite grids of any finite dimension $d$.
A natural way to define the odd independence density $\rhodd$
 in higher-dimensional grids is to take ${\displaystyle 
 \sup_S \limsup_{H_n\,, \ n\to\infty} } |S\cap V(H_n)|/n^d$ where all \sois s $S$
 of the grid are considered, and the current $H_n$ ranges over all
 $d$-dimensional $P_n\Box\cdots\Boxx P_n$ in the grid.

Let us consider first the chromatic number with respect to distance 2.

\bpn
For the $d$-dimensional infinite grid $\gdinf:=P_\infty \Boxx \cdots \Box
 P_\infty$ on the vertex set $\zzz^d$ we have $
 \chi((\gdinf)^2) = 2d+1$.
\epn

\bpf
It is clear that $\chi((\gdinf)^2) \geq \omega((\gdinf)^2) \geq
 \Delta(\gdinf)+1 = 2d+1$.
Consider now the coloring $\vp : \zzz^d \to \{ 0,1,\dots,2d \}$
 on $\gdinf$ determined by the rule $\vp (a_1,\dots,a_d) = \bigl( \,
  \sum_{i=1}^d i\cdot a_i \bigr)$ (mod $2d+1$).
This assignment is a proper vertex coloring of $(\gdinf)^2$,
 because the weighted sums defined above in the neighborhood of
 any vertex $v$ differ by at most $2d$, taking all values $\vp(v)-d,
 \vp(v)-d+1, \dots, \vp(v)-1, \vp(v)+1, \vp(v)+2, \dots, \vp(v)+d$
  modulo $2d+1$.
Thus, the assertion follows.
\epf

Since the inequality $\soc(G) \leq \chi(G^2)$ holds for
 all graphs $G$, it follows that
 $\soc(\gdinf) \leq 2d+1$ and then $\rhodd(\gdinf) \geq 1/(2d+1)$.
It turns out that these estimates are very far from being tight,
 as $\soc$ does not depend on $d$.

\btm
\label{t:grid-d}
We have $\soc(\gdinf) = 3$ for all $d\geq 1$.
As a consequence, $\rhodd(\gdinf) \geq 1/3$.
\etm

\bpf
Assume that $\vp$ is a \sokk\ on $\gdinf$.
The vertex degrees are equal to $2d$, hence
 due to the parity conditions each vertex $v$
 has two neighbors $x,y$ with $\vp(x)\neq \vp(y)$.
Since $\vp$ is also a proper vertex coloring,
 $\vp(v)$ is a third color, implying $\soc \geq 3$.
The proof will be done if we show that three colors are sufficient.

Let $v=(a_1 , \dots , a_d)$ be any vertex of $\gdinf$.
If $d$ is odd, we define the color of $v$ as
 $\vp(v) := (a_1 + \dots + a_d)$ [mod 3]\,.
Then, modulo 3, $v$ has $d$ neighbors of color $\vp(v)+1$
 and also $d$ neighbors of color $\vp(v)-1$, hence
 both multiplicities are odd, and a
 \sokk\ of $\gdinf$ with three colors is obtained.

\begin{figure}[!ht]
\centering
\begin{subfigure}{0.45\textwidth}
\centering
\includegraphics[width = \textwidth]{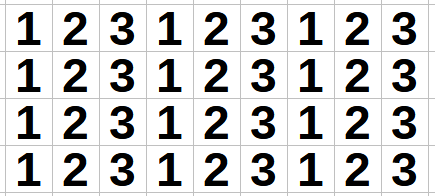}
\caption{}
\label{fig:h1}
\end{subfigure}\hfill
\begin{subfigure}{0.45\textwidth}
\centering
\includegraphics[width = \textwidth]{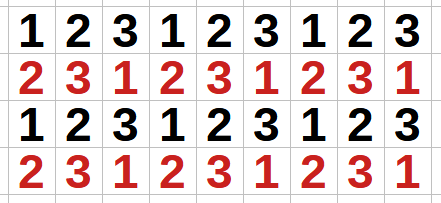}
\caption{}
\label{fig:h2}
\end{subfigure}
\caption{Strong odd 3-coloring of $P_\infty\Boxx P_\infty$, cellular representation. (a) Odd coloring in the first coordinate, (b)~modified coloring in levels of odd index.}
\label{fig:grid-3coloring}
\end{figure}

Assume from now on that $d$ is even.
We introduce two auxiliary functions:

 $$
   h'(v) := ( a_1 + \dots + a_{d-1} ) \ [\mathrm{mod} \ 3]\,, \qquad
   h''(v) := a_d \ [\mathrm{mod} \ 2] \,;
 $$
  and then define

 $$
   \vp(v) := ( \, h'(v) + h''(v) \, ) \ [\mathrm{mod} \ 3] \,.
 $$
Purely $h'$ assigns color $h'(v)-1$ to $d-1$ neighbors of $v$,
 also color $h'(v)+1$ to $d-1$ neighbors of $v$
 (both multiplicities are odd),
 and color $h'(v)$ to just two neighbors of $v$ (those two
 which agree with $v$ in all of the first $d-1$ coordinates),
An illustration for $d=2$ is exhibited in Figure \ref{fig:grid-3coloring};
 the cells modified with $h''$ are indicated with red numbers.

If $a_n$ is even, then modifying $h'$ to $\vp$
 keeps the color of $v$ unchanged, as well as the colors of
 those vertices whose last coordinate equals $a_d$.
This includes $2d-2$ neighbors of $v$ in the
 colors $h'(v)-1$ and $h'(v)+1$.
The last two neighbors of $v$ change their color
 from $h'(v)$ to $h'(v)+1$.
As a result, $d+1$ neighbors of $v$ (an odd number)
 have color $\vp(v)+1$, and $d-1$ neighbors
 (also odd) have color $\vp(v)-1$.

On the other hand,
 if $a_n$ is odd, then modifying $h'$ to $\vp$
 changes the color of $v$ from $h'(v)$ to $h'(v)+1$
 (computed modulo 3, of course), as well as the colors of
 those vertices whose last coordinate equals $a_d$.
This includes $2d-2$ neighbors of $v$ in the
 colors $h'(v)-1$ and $h'(v)+1$.
The last two neighbors of $v$ change their color
 from $h'(v)$ to $h'(v)+1$.
As a result, $d+1$ neighbors of $v$ (an odd number)
 have color $\vp(v)+1$, and $d-1$ neighbors
 (also odd) have color $\vp(v)-1$.
Further, all vertices $w$ having last coordinate $a_n$
 undergo the same modification.
As a consequence, among them $d-1$ have
 $\vp(w)=\vp(v)+1$ and the same number of them have $\vp(w)=\vp(v)-1$.
Since the last two neighbors of $v$ have an even last coordinate,
 they keep their color $h'(v)=\vp(v)-1$.
Thus, in this case $d+1$ neighbors of $v$
 have color $\vp(v)-1$, and $d-1$ neighbors
 have color $\vp(v)+1$.
\epf

\subsection{Planar grids}

Here we briefly consider the function $f(n,k):=\soc(P_k\Boxx P_n)$,
 where $k$ is fixed and $n$ is arbitrary.
Of course, $k=1$ just means the graph $P_n$, already settled
 as $\soc(P_n)=\chi((P_n)^2)=3$ for all $n\geq 3$, with $\sid(P_n)=\ceil{n/3}$.

  \btm 
\label{p:stripe}
The following holds for the function $f(n,k)$ defined above.
  \begin{itemize}
  \item[$(i)$] $f(n,k)\leq 5$ for all $n,k$;
  \item[$(ii)$] $f(n,k)\geq 3$ for all $n,k\geq 2$;
  \item[$(iii)$] $f(n,k)\geq 4$ for all $n,k\geq 4$;
  \item[$(iv)$] $f(n,k)\leq 4$ whenever $2\mid k$ and $k-1 \mid n$;
  \item[$(v)$] $f(k,k)=5$ for all odd $k\geq3$.
  \end{itemize}
\etm

\bpf
Throughout we use a cellular representation of $G=P_n\Boxx P_k$. Consider the following vertex coloring $\varphi:V(G)\to\{0,1,2,3,4\}$. For any vertex $v=(a,b)$ of $G$ let us define $\varphi(v)=2\cdot (a\Mod{5})+(b\Mod{5})$, addition modulo 5. One readily observes that if $\varphi(v)=x$ then its neighbors are assigned with distinct values from the set $\{x-2,x-1,x+1,x+2\}$, addition and subtraction modulo 5. Thus $\varphi$ is a proper coloring of $G^2$, which in turn gives that $\varphi$ is a strong odd $5$-coloring of $G$. This implies $(i)$.

  Concerning $(ii)$, two colors are not sufficient for any pair $(n,k)$ with $n,k\geq 2$ because a corner vertex and its two neighbors have to receive mutually distinct colors under any strong odd coloring.

  Let $n,k\geq 4$, and suppose that $G$ admits a strong odd coloring with color set $\{1,2,3\}$. Upon permuting colors, we may assume that $f(1,1)=f(2,2)=1$, $f(1,2)=2$ and $f(2,1)=3$. Consequently, $f(1,3)=f(3,1)=1$ and $f(2,3)=f(1,4)=2$. However, then it must be that $f(3,2)=2$ as well, which in turn leaves us without an option for $f(4,1)$. The obtained contradiction proves $(iii)$.

\begin{figure}[ht]
  \centering
		\includegraphics[scale=0.35]{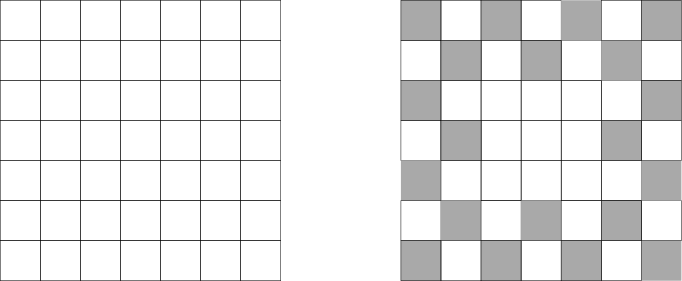}
	\caption{A $7\times 7$ square table (left), and its frame (right).}
	\label{fig:frame}
\end{figure}

  Concerning $(iv)$, let us first construct a strong odd $4$-coloring of $P_{2t-1}\Boxx P_{2t}$.  We introduce the following ad-hoc terminology. A \textit{peripheral cell} of the table is any $1\times 1$ cell $(a,b)$ such that $a\in\{1,2,2t-2,2t-1\}$ or $b\in\{1,2,2t-1,2t\}$ (or both). Note that the non-peripheral cells comprise a `central' $(2t-5)\times (2t-4)$ (sub)table. For any $(2s-1)\times (2s-1)$ square (sub)table within the $(2t-1)\times (2t)$ rectangle, its \textit{frame} consists of the cells depicted as grey in Figure~\ref{fig:frame}.

\begin{figure}[ht]
  \centering
		\includegraphics[scale=0.35]{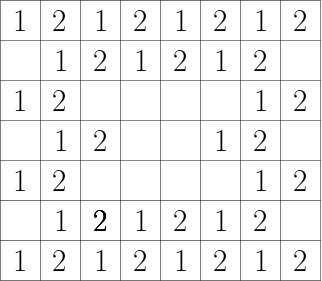}
	\caption{A partially colored $7\times 8$ rectangular table. Apart from the cells in even-numbered rows from the first and last column, every other peripheral cell is assigned with a color from the set $\{1,2\}$. Additionally, some non-peripheral cells are assigned with a color from the set $\{1,2\}$ --- these cells are peripheral (located in even-numbered rows from the first and last column) in regard to the central $3\times 4$ rectangular (sub)table.}
	\label{fig:coloredframe}
\end{figure}

  Now look at the $(2t-1)\times (2t)$ rectangle (in the cellular representation of $P_{2t-1}\Boxx P_{2t}$) as if comprised of two partially overlapping $(2t-1)\times (2t-1)$ square tables: a left one consisting of the cells $(a,b)$ with $1\leq b\leq 2t-1$, and a right one consisting of the cells $(a,b)$ with $2\leq b\leq 2t$. Assign to each frame cell of the left $(2t-1)\times (2t-1)$ square the color $1$. Similarly, assign to each frame cell of the right $(2t-1)\times (2t-1)$ square the color $2$ (cf.\ Figure~\ref{fig:coloredframe}).

\begin{figure}[ht]
  \centering
		\includegraphics[scale=0.27]{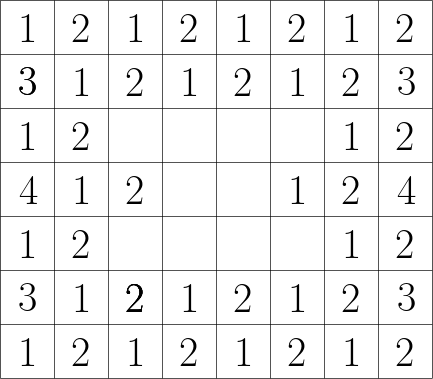}
	\caption{A partially colored $7\times 8$ rectangular table. Every frame cell is assigned with a color from the set $\{1,2\}$. The remaining peripheral cells are assigned with a color from the set $\{3,4\}$. Additionally, the cells that are in even-numbered rows from the first and last column in regard to the central $3\times 4$ rectangular (sub)table are assigned with a color from the set $\{1,2\}$. }
	\label{fig:fullycoloredframe}
\end{figure}

We assign to the remaining uncolored peripheral cells a color from the set $\{3,4\}$ as follows: for the uncolored cells in the first column (the cells $(a,b)$ with $a\in\{2,4,2t-2\}$ and $b=1$) we alternate between the color $3$ and the color $4$ as the row coordinate of the cell increases; we color analogously the uncolored cells in the last column (the cells $(a,b)$ with $a\in\{2,4,2t-2\}$ and $b=2t$) (cf.\ Figure~\ref{fig:fullycoloredframe}).

Note that the uncolored cells comprise the central $(2t-5)\times (2t-4)$ (sub)table, excepting the peripheral cells of this subtable which come from its even-numbered row in its first/last column. We iterate the `frame coloring' procedure with permuted color: instead of colors $1$ and $2$ we use $3$ and $4$, respectively, and vice versa. In more detail, we look at this central $(2t-5)\times (2t-4)$  rectangle (its peripheral cells in even-numbered rows from the first/last column are already colored with $1$ or $2$) as comprised of two partially overlapping $(2t-5)\times (2t-5)$ squares, and color their respective frames using colors $3$ and $4$ (the color $3$ for the frame of the left square, and the color $4$ for the frame of the right square). And so on, until every cell gets colored. The final outcome is a strong odd coloring of $P_{2t-1}\Boxx P_{2t}$ such that every cell in the first and last row is assigned with a color from the set $\{1,2\}$.

 The mentioned extra feature for the colors of the first/last row enables us to generate a strong odd $4$-coloring for any $P_{m(2t-1)}\Boxx P_{2t}$. Indeed, given $P_{m(2t-1)}\Boxx P_{2t}$, decompose it into $m$ copies of $P_{2t-1}\Boxx P_{2t}$, and enumerate these copies as first, second, $\ldots,$ $m$-th so that the last row of the $(i-1)$-st copy of $P_{m(2t-1)}\Boxx P_{2t}$ is adjacent to the first row of the $i$-th copy for any $i=1,2,\ldots,m$. For the first copy of $P_{m(2t-1)}\Boxx P_{2t}$ use the already constructed strong odd $4$-coloring. Once the coloring of the $(i-1)$-st copy has been completed, color the $i$-th copy of $P_{2t-1}\Boxx P_{2t}$ by switching the roles played by  colors $1$ and $3$, as well as the roles played by colors $2$ and $4$ in the coloring of the previous copy of $P_{2t-1}\Boxx P_{2t}$. This gives a strong odd $4$-coloring of $P_{m(2t-1)}\Boxx P_{2t}$, which proves $(iv)$.

Finally, let us show that $P_{2t+1}\Boxx P_{2t+1}$ requires five colors for any strong odd coloring. Consider the pairs of cells that cross the main diagonal: the first such pair is comprised of the cells $(1,2)$ and $(2,1)$, the second pair of the cells $(2,3)$ and $(3,2)$, and so on, the last pair is comprised of the cells $(2t,2t+1)$ and $(2t+1,2t)$. The total number of considered pairs is even (namely, $2t$). We argue by contradiction, and suppose a strong odd $4$-coloring is possible. Then the two cells comprising the first pair must receive distinct colors (in view of the corner cell $(1,1)$). Consequently, the two cells forming the second pair must be colored the same (in view of the cell $(2,2)$). Consequently, the two cells forming the third pair must receive distinct colors (in view of the cell $(3,3)$), and so on. We conclude that the two cells forming the last considered pair must receive the same color, since the total number of considered pairs of cells is even, $2t$, as already noted. However, this cannot be the case because then the corner cell $(2t+1,2t+1)$ would have its only two neighboring cells colored the same.
    \qed

  \bigskip

  Note that $(iii)$ is tight as $f(3,4)=3$. In view of Theorem \ref{p:stripe}, for all $n,k\geq4$ it holds that $f(n,k)$ equals $4$ or $5$, and either value is attained for infinitely many pairs $(n,k)$. It is our belief that $f(k,k)=5$ for all even $k\geq4$; however, we are unable to provide a succinct argument for this.
\epf

\section{Chessboard graphs}
\label{s:chess}

The chessboard graphs are defined on the vertex set $\{1,\dots,n\} \times \{1,\dots,n\}$, where adjacencies correspond to the moving rules of chess.
We call $n$ the size of the graph in question, referring to the side length of the board on which it is taken.
Such graphs are defined for the king, queen, rook, bishop, and knight.
In some of them the determination of $\sid$ is easy, in some others it leads to a challenging open problem.

As in the preceding section, also here a cellular representation can be applied.
It is helpful in visualization; equivalence between the
 $n^2$ vertices and the $n^2$ cells (as well as in the infinite case) is established by using coordinate pairs $(i,j)$.

\paragraph{King.}

The King graph  is the same as the strong product $P_n \strg P_n$.

\bpn
\label{p:king}
If $G$ is the King graph of size $n$, then $\sid(G)=\ceil{n/3}^2$, and $\soc(G)=9$ for all $n\geq 3$.
Moreover, the largest odd independent set is unique if and only if $n\equiv 1 \ (\mathrm{mod}~3)$.
\epn

\bpf
Recall that every independent set in $G^2$ is an \sois, in any graph $G$.
We first show that in the King graph the converse implication is also valid: if $S$ is an \sois\ in the King graph $G$, then $S$ is also an independent set in $G^2$.
Once proved, we can infer $\sid(G) = \alpha(G^2)$.


The eight neighbors of a vertex $(a,b)$ are of the form $(a\pm1,b\pm1)$.
Among them, $S$ can contain at most one $(a',b')$ with $|a-a'|+|b-b'|=1$, because e.g.\ $(a-1,b)$ forms a forbidden pair with $(a+1,b)$, and is adjacent to both $(a,b-1)$ and $(a,b+1)$.
From the other four with $|a-a'|+|b-b'|=2$, the only non-forbidden pairs are the opposite diagonal ones.
However, if $(a-1,b-1), (a+1,b+1)\in S$, then this vertex pair dominates every neighbor of $(a,b)$ with $|a-a'|+|b-b'|=1$, thus $|S\cap N[(a,b)]|\leq 2$.
By the parity condition on $S$, this implies that every vertex has at most one neighbor in $S$.

Independence in $G^2$ means that the sets $N[v]$ are mutually disjoint for all $v\in S$.
Then for any two $(a,b),(a',b')\in S$ we have $|a-a'|\geq3$ or $|b-b'|\geq3$ or both.
Viewing $S$ as the vertex set of an auxiliary $K_{|S|}$, orient a blue arc from $(a,b)$ to $(a',b')$ if $a'\geq a+3$; and orient a red arc if $|a-a'|\leq 2$ and $b'\geq b+3$.
Clearly, no monochromatic directed path can have more than $\ceil{n/3}$ vertices, and $\ceil{n/3}$ is attained for $n\equiv 1 \ (\mathrm{mod}~3)$ by the unique sequence $1,4,7,\dots,n$.
If $b_v,r_v$ are the lengths of the longest blue and red path starting from $v$, then for any two $v,w\in S$ we have $b_v\neq b_w$ or $r_v\neq r_w$ or both (due to the fact that there is a blue or red arc between $v$ and $w$).
Thus, $|S|\leq \ceil{n/3}^2$, as claimed.
Also, the assertion concerning uniqueness follows.

A \sokk\ with 9 colors is obtained by the color assignment $(a,b)\mapsto (a$ [mod 3]) + $3\cdot (b$ [mod 3]).
Fewer colors are not enough because, as seen above, already a $3\times 3$ subsquare requires mutually distinct colors on its nine vertices.
\epf

The same argument yields $\sid(P_p \strg P_q)=\ceil{p/3}\cdot\ceil{q/3}$ for all $p,q\geq 1$, with a unique largest \sois\ if and only if both $p$ and $q$ are congruent to 1 modulo 3.
Also the \sokk\ with nine colors, as defined in the proof, works for all $p,q\geq 3$.

The above results can be generalized as follows.
The $r$-King graph ($r\in\nnn$) of size $n$ has the same vertex set of $n^2$ vertices as the other chessboard graphs.
Two vertices $(a,b),(a',b')\in\{1,\dots,n\}\times\{1,\dots,n\}$ are adjacent if and only if $|a-a'|\leq r$ and $|b-b'|\leq r$.
The King graph represents the particular case $r=1$.

\btm
The $r$-King graph of any size $n$ has $\sid=\ceil{\frac{n}{2r+1}}^2$.
Moreover, $\soc=(2r+1)^2$ if $n\geq 2r+1$, and $\soc=n^2$ if $n\leq 2r$.
\etm

\bpf
Let $G$ be the $r$-King graph of size $n$.
Recall that in the $n\times n$ square board representation the rows and columns are numbered from 1 to $n$,
 hence vertex $(a,b)$ of $G$ naturally corresponds to the cell in the intersection of column $a$ and row $b$,
 making vertices and cells equivalent.

\begin{figure}[!ht]
\begin{center}
\resizebox{8cm}{!}{%
\includegraphics
{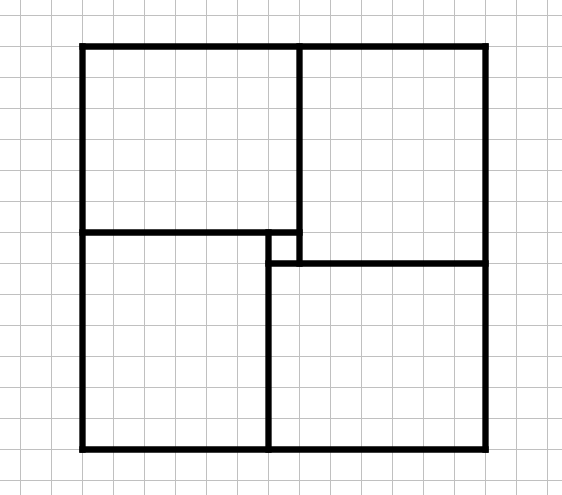}
  }
    \caption{Partition of a $(2r+1)\times(2r+1)$ subsquare ($r=6$).}
    \label{fig:2r+1-split}
\end{center}
\end{figure}

\begin{figure}[!ht]
\begin{center}
\resizebox{12cm}{!}{%
\includegraphics
{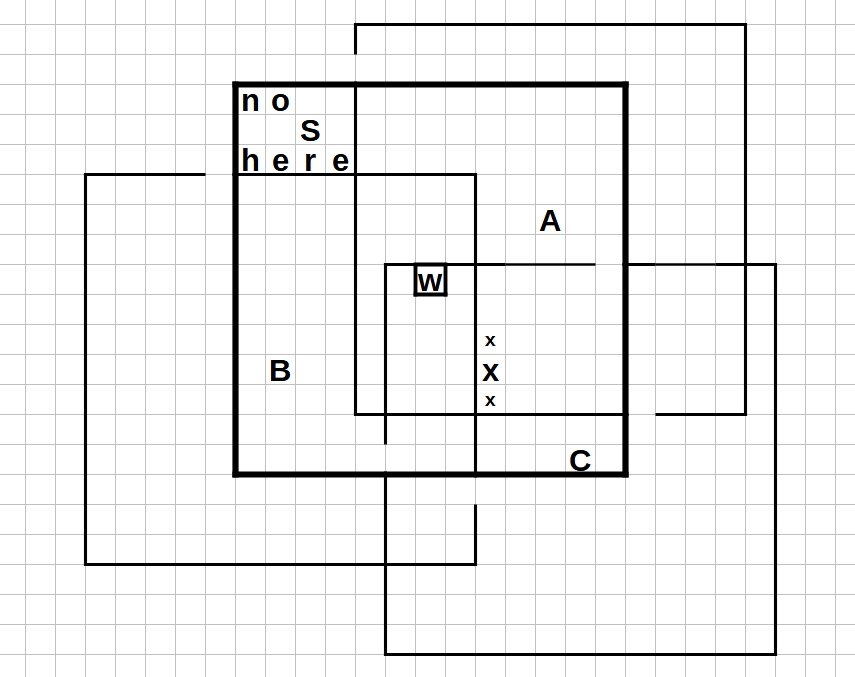}
  }
    \caption{Thick lines: border of the $(2r+1)\times(2r+1)$ subsquare with center $w$; lines of half thickness: boundaries of the neighborhoods of $A,B,C$; some cells next to $N[B]$ are indicated with $X$ or ${\scriptstyle X}$.}
    \label{fig:2r+1-ABC}
\end{center}
\end{figure}

Consider any \sois\ $S$ in $G$.
We prove that no $(2r+1)\times(2r+1)$ subsquare can contain more than one vertex from $S$.
Suppose for a contradiction that a $(2r+1)\times(2r+1)$ subsquare $Q$ contains more than one member of $S$.
The central cell, say $w$, of $Q$ is not in $S$, because it is not independent of any cell of $Q$.
We split $Q-w$ into four $r\times(r+1)$ rectangles as shown in Figure \ref{fig:2r+1-split}.
Inside each rectangle any two cells are dependent, hence $S$ has at most one vertex there.
Further, exactly one of those four rectangles must be free of $S$, otherwise $w$ would be dominated by two or four members of $S$, contrary to odd independence.
We assume without loss of generality that the upper left rectangle is blank.
Let the three vertices $A,B,C\in S$ adjacent to $w$ be arranged as indicated in Figure \ref{fig:2r+1-ABC}; that is, $A$, $B$, $C$ are placed in the upper right, lower left, and lower right rectangle, respectively.

We use the notation $<_h$, $\leq_h$ to indicate relative horizontal positions with respect to columns; for instance, $B<_h w\leq_h C$ is implied by the chosen rectangle partition.
We write $<_v$ and $\leq_v$ with a similar meaing in the vertical direction; for instance $C<_v A$.

Since $A$ and $B$ are independent in $G$, the horizontal or vertical distance between $A$ and $B$ must exceed $r$.
We assume that the former is valid (or both), with the above notation we can write this relation as $N[B]<_hA$.

We are now most interested in the area $R=(N[A]\cap N[C])\smin N[B]$.
The cells in $R$ are dominated by $A$ and $C$, but not by $B$.
We are going to prove that this $R$, more precisely its column next to the right edge of $N[B]$, contains a cell $X$ such that $N[X]\subseteq N[A]\cup N[C]$.
If this property holds, then $A,C$ is a forcing pair for $X$, which contradicts the assumption that $S$ is an \sois, due to Lemma \ref{l:force}.
Thus, $|S\cap Q|\leq 1$ will then be proved.

Assuming that this property is valid, it immediately implies that in a \sokk\ all vertices must have distinct colors if $n\leq 2r+1$, and also at least $(2r+1)^2$ colors are needed if $n\geq 2r+1$.
On the other hand, $(2r+1)^2$ colors are always enough, as shown by the color assignment $(a,b) \mapsto (a$ mod $r) + r\cdot (b$ mod $r)$.

For the rest of the proof, concerning $\sid$, we apply an argument analogous to the second part of the proof of Proposition \ref{p:king}.
We define the oriented blue-red edge-colored complete graph $K_{|S|}$, orienting a blue
arc from $(a, b)$ to $(a', b')$ if $a' \geq a + 2r+1$, and a red arc if $|a - a'| \leq 2r$
and $b' \geq b + 2r+1$. Then no monochromatic directed path can have more
than $\ceil{\frac{n}{2r+1}}$ vertices, and 
for $n \equiv 1$ (mod $2r+1$) this bound is attained by the unique
sequence 1, $2r+2$, $4r+3, \dots , n$.
Also here, if $b_v, r_v$ are the lengths of the longest blue and red
path starting from vertex $v$, then for any two $v,w \in S$ we have $b_v \neq b_w$ or $r_v \neq r_w$ or both.
Consequently, $|S| \leq \ceil{\frac{n}{2r+1}}^2$.
\epf

\paragraph{Rook.}

The Rook graph of size $n$ is $K_n \Boxx K_n$.

\bpn
$\sid(K_n \Boxx K_n)=1$ for every $n\geq 1$; and $\soc(K_n \Boxx K_n)=n^2$.
\epn

\bpf
The open neighborhood of every vertex is isomorphic to $2K_{n-1}$, hence the graph is claw-free, therefore $\sid(G)=\alpha(G^2)$, as proved in \cite{part-1}.
Since $G$ has diameter 2, we obtain $\alpha(G^2)=1$.
Then also $\soc(K_n \Boxx K_n)=n^2$ follows immediately.
\epf

One can also observe the extensions $\sid(K_p \Boxx K_q)=1$ and $\soc(K_p \Boxx K_q)=pq$ for all $p,q\geq 1$.

As for the king, also for the rook, finite and infinite
 $r$-Rook graphs can be defined for any $r\in\nnn$.
In this case the horizontal or vertical move of the rook
 is limited to distance $r$.
Taking the representation with coordinate pairs,
 this means adjacency of $(a,b)$ and $(a',b')$ if and only if
 either $a=a'$ and $|b-b'|\leq r$,
 or $b=b'$ and $|a-a'|\leq r$.
For the infinite $r$-Rook graph we have general estimates
 in which the ratio of the upper and lower bound tends to 3/2 as $r$ gets large.

\btm
\label{t:r-rook}
The odd density $\rhodd$ of the infinite $r$-Rook graph satisfies

 $$
   \frac{1}{2r+1} \leq \rhodd \leq \frac{3}{4r+3} \,.
 $$
\etm

\begin{figure}[!ht]
\begin{center}
\resizebox{8cm}{!}{%
\includegraphics
{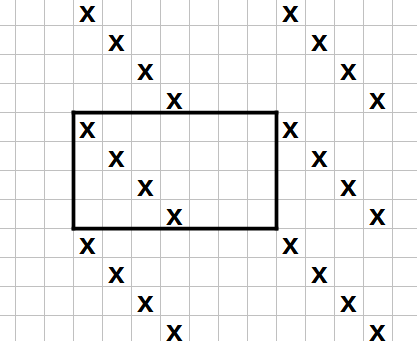}
  }
    \caption{Detail of the infinite 3-Rook graph; {\large\sf x} = vertex selected for $S$.}
    \label{fig:r-rook}
\end{center}
\end{figure}

\bpf
The upper bound is a consequence of Proposition \ref{p:reg-clawfree-asymp},
 taking into account that the graph is $K_{1,5}$-free and
 has maximum degree $4r$.

For the lower bound we arrange $r+1$ rooks in a diagonal position
 and leave the next $r$ columns blank.
We repeat this pattern vertically modulo $r+1$, and horizontally modulo $2r+1$.
Illustration for $r=3$ is given in Figure~\ref{fig:r-rook}.
From any $(r+1)\times (2r+1)$ rectangle, exactly $r+1$ cells (vertices)
 are selected, hence the density equals $\frac{1}{2r+1}$.
Each cell in the blank columns is dominated just once from $S$,
 and the blank cells in the other columns have precisely three neighbors in $S$.
Thus, $S$ is an \sois .
\epf

\paragraph{Bishop.}

In the Bishop graph of any size, two vertices $(a,b)$ and $(a',b')$ are adjacent if and only if $a-a'=b-b'$ or $a-a'=b'-b$.
It has two components, the connected white and black graphs, which are isomorphic if and only if $n$ is even.
Moreover, each component (black or white) of the infinite Bishop graph
 is isomorphic to the infinite Rook graph.

\bpn
\label{p:bishop}
The Bishop graph of any size $n\geq 2$ has $\sid=2$ and $\soc=\ceil{n^2\!/2}$.
\epn

\bpf
In each component, the open neighborhood of each vertex is the disjoint union of two complete graphs (or is just one complete graph if the vertex is in a corner), hence the graph is claw-free.
Consequently $\sid
(G)=\alpha(G^2)$, which is now equal to 2 because both components have diameter 2 (or 1 if $n=2$).

The larger component has $\ceil{n^2\!/2}$ vertices, each of which is a maximal \sois\ and so requires a private color in its component.
Assigning this way, and repeating those colors (except one if $n$ is odd) in the other component, a \sokk\ is obtained.
\epf

For the bishop, the above result does not extend directly to a $p\times q$ board, because the diameter of the corresponding graph can be large, depending on the ratio between $p$ and $q$.
This admits larger odd independence, and hence also \sokk\ with fewer colors.

Similarly to the case of rook, in the $r$-Bishop graphs ($r\in\nnn$) the (diagonal) move of a bishop is limited to distance at most $r$.
Expressing this condition with coordinates, adjacency of $(a,b)$ and $(a',b')$ is taken under the restriction $|a-a'|\leq r$, or equivalently $|b-b'|\leq r$.
In this way the infinite $r$-Bishop graph has two components, each isomorphic to the infinite $r$-Rook graph.

\btm
The odd density $\rhodd$ of the infinite $r$-Bishop graph satisfies

 $$
   \frac{1}{2r+1} \leq \rhodd \leq \frac{3}{4r+3} \,.
 $$
\etm

\begin{figure}[!ht]
\begin{center}
\resizebox{\textwidth}{!}{%
\includegraphics
{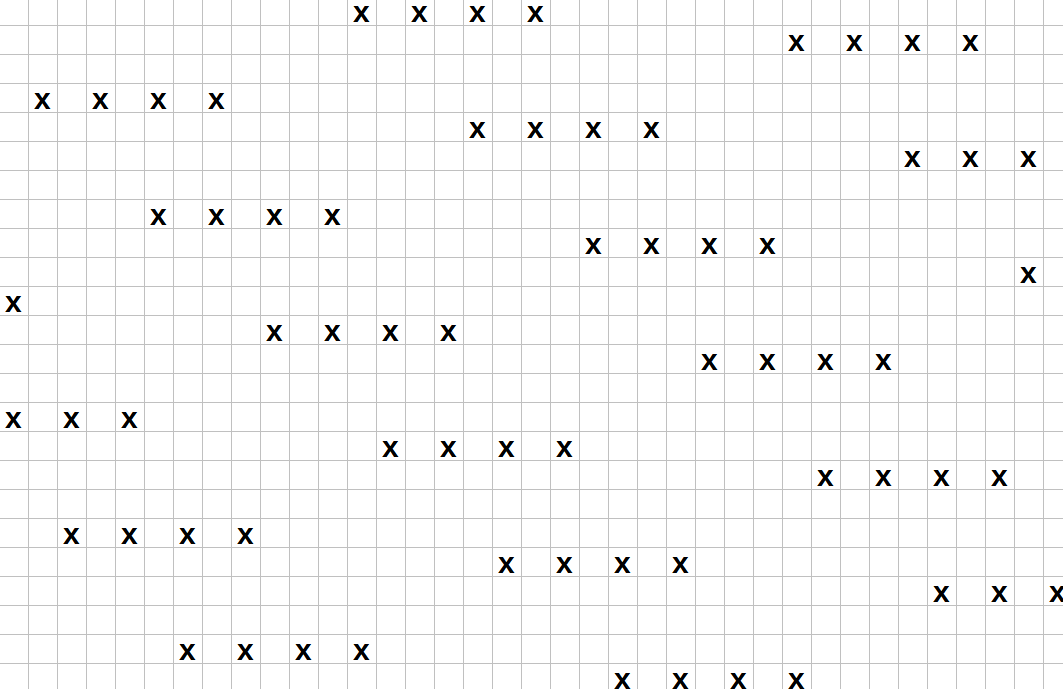}
  }
    \caption{Odd-independent set of density 1/7 in one component of the infinite 3-Bishop graph.}
    \label{fig:r-bishop}
\end{center}
\end{figure}

\bpf
Since each connected component (black or white) is isomorphic to the $r$-Rook graph,
 the inequalities follow by Theorem \ref{t:r-rook}.
Illustration for $r=3$ is given in Figure \ref{fig:r-bishop}.
The easiest way to check density is to observe that each up-diagonal line
 of the component in question contains exactly six gaps
 between any two selected cells.
An analogous selection can be done in the other component, too;
 but there are several ways to combine the two selections.
\epf

\paragraph{Queen.}

The edge set of the Queen graph is the union of those of the Rook and Bishop graphs.
It has diameter 2, and every vertex is in the intersection of at most four cliques, hence it is $K_{1,5}$-free.

\bpn
\label{p:queen-10}
Every Queen graph of size $n\leq 10$ has $\sid=1$, and consequently also $\soc=n^2$.
\epn

\bpf
This has been checked by running a computer code.
\epf

\bpn
\label{p:queen-inf}
The infinite Queen graph has either $\sid=\infty$ or $\sid=1$ but
 no other finite value is possible.
\epn

\bpf
Suppose for a contradiction that $\sid=k$, where
 $1<k<\infty$, and consider an \sois\ $S$ of $k$ vertices.
Let $(a_1,b_1), (a_2,b_2) \in S$ be the two vertices such that
 $a_1$ is the largest first coordinate and
 $|b_1-b_2|$ is the maximum difference from $b_1$ in $S$.
We consider the vertex $v=(a,b)$ where $b=b_2$ and
  \begin{itemize}
   \item $(a_1+b_1+b_2,b_2)$ if $b_2>b_1$, hence $b-a = b_1-a_1$;
   \item $a=a_1+b_1-b_2$ if $b_2<b_1$, hence $a+b = a_1+b_1$.
  \end{itemize}
This $v$ has precisely two neighbors in $S$, namely $(a_2,b_2)$ horizontally
 and $(a_1,b_1)$ diagonally, because the vertical line and the
 other diagonal passing through $v$ is disjoint from $S$.
This contradicts the assumption that $S$ is an \sois.
\epf

Clearly, the infinite Queen graph has $\soc = \infty$, but this does not allow any conclusion concerning its $\sid$ whether it is 1 or $\infty$.

\paragraph{Knight.}

A general vertex (at distance at least 2 from the edge) of the Knight graph, say $(a,b)$, has eight neighbors:

\begin{center}
$(a-2,b+1)$, $(a-2,b-1)$, $(a+2,b+1)$, $(a+2,b-1)$, 

$(a-1,b+2)$, $(a-1,b-2)$, $(a+1,b+2)$, $(a+1,b-2)$.
\end{center}

\bpn
\label{p:knight}
As $n\to\infty$, the Knight graph of size $n$ has $\frac{7}{16}\,n^2 - O(n) \leq \sid \leq \frac{7}{15}\,n^2 + O(n)$.
\epn

\bpf
The asymptotic upper bound follows as the particular case $\Delta=8$ of Proposition \ref{p:reg-clawfree-asymp} $(i)$.
The error term can be obtained by observing that the vertices of degree less than 8 are located in the first and last two rows and columns, hence $z_j = 4\cdot (2n-4) = O(n)$
using the notation introduced in the proof of that proposition. 

\begin{figure}[!ht]
\begin{center}
\resizebox{\columnwidth}{!}{%
\includegraphics
{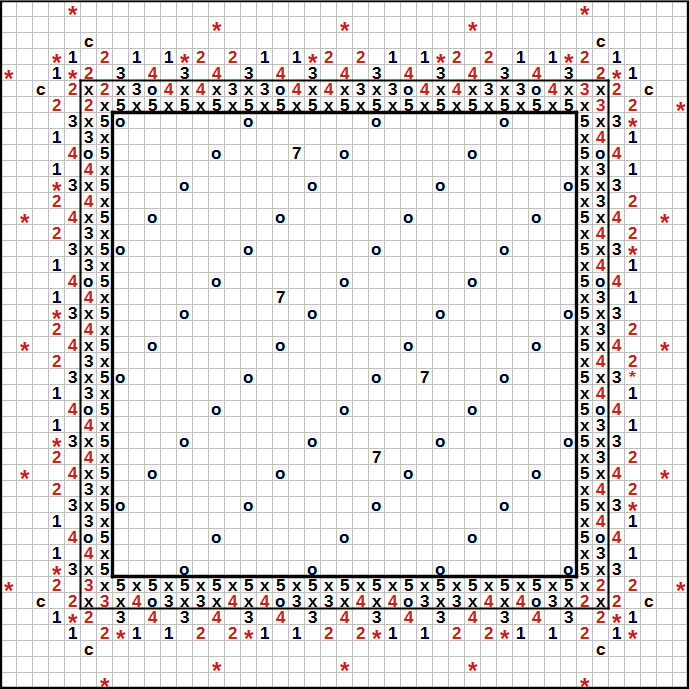}
  }
    \caption{Knight graph, odd independent density $7/16-o(1)$.}
    \label{fig:knight}
\end{center}
\end{figure}

We prove the lower bound by construction; an illustrative example with $n=43$ is exhibited in Figure \ref{fig:knight}, where again the cells of the board represent the $n^2$ vertices of the graph.
Recall that the number of both the black and white cells is essentially the same, namely $\floor{n^2\!/2}$ and $\ceil{n^2\!/2}$.

The construction begins with selecting approximately $7/8$ of black vertices.
The pattern of non-selected cells is indicated with {\large\textsf{\textbf{o}}} in the big central area.
In rows of odd index all black cells are selected, in rows of even index one black is skipped after three black ones; and this is shifted diagonally as shown in the figure.
As a result, every white cell of the central area is knight-connected to precisely one non-selected black cell, hence dominated by exactly seven selected cells.
Four white cells labeled 7 intend to express this fact.
One can observe that from each of them the closest {\large\textsf{\textbf{o}}}\,---and only that one---is reached by an L-shape move of the knight.

The second step is to remove all selections from the first and last five rows and columns.
In rows/columns of indices 6, 7, $n-6$, $n-5$ the selections are kept and are indicated with {\large\textsf{\textbf{x}}} in the figure.
The numbers in rows/columns from 4 to 7 and from $n-6$ to $n-3$ indicate how many selected black cells dominate the white cell in question.
Many of them are odd (1, 3, 5) but some are even (2, 4), written in red.

The third step is to select further black cells in the rows and columns of indices from 2 to 5 and from $n-4$ to $n-1$, as indicated by the {\Large $\star$} signs.
Each red even number is dominated by ecactly one newly selected cell, hence its domination is increased from 2 to 3 or from 4 to 5.
In addition, near some cells, some 3-labeled cells get dominated by precisely two newly selected cells, hence their domination is increased from 3 to 5; these labels 3 are written in red.

As a side effect, two cells near each corner (in row/column 3 and $n-2$) were not dominated before, but now have exactly two newly selected neighbors in the Knight graph.
These are indicated with {\large\textsf{\textbf{c}}}.
In the last (fourth) step we give them just one further domination in row/column 1 or $n$, indicated by the {\Large $\star$} sign, increasing its selected neighbors to 3.

In this way an \sois\ of cardinality at least $\frac{7}{16}(n-10)^2 = \frac{7}{16}\cdot n^2 - O(n)$ has been selected, yielding the claimed lower bound.
\epf

\section{Other types of grids}
\label{s:tri-hex}

The two propositions of this short section open a beautiful area
 on the border of discrete geometry and graph theory, namely
 independence and odd independent density in tilings of the plane,
 and even more generally on surfaces.
We now determine the odd independence ratio
 in the infinite triangular and hexagonal grids.
As above, also here the odd independence ratio is denoted by $\rhodd$.

\bpn
The infinite triangular grid has $\rhodd = 1/3$ and $\soc = 3$.
\epn

\bpf
The triangular grid is a 6-regular graph whose chromatic number is 3.
Take a proper 3-coloring $\vp$, and let $S$ be any one of its color classes.
Consider any $v\notin S$.
Then $N(v)$ induces a 6-cycle, which is properly colored with the two colors distinct from $\vp(v)$, one of them being the color of $S$.
Hence, $v$ has precisely three neighbors in $S$, therefore $S$ is an \sois\ of density 1/3, and $\vp$ is a \sokk.

Denser \sois s---and also denser independent sets---cannot exist, because the grid obviously admits a vertex partition into triangles (e.g., the grid can be split into parallel stripes each of which induces a subgraph isomorphic to $(P_\infty)^2$), and any independent set can contain at most one vertex from each triangle.
\epf

We note that the maximum density of odd subgraphs is 6/7, hence much larger than $\rhodd$.

\bpn
The infinite hexagonal grid has $\rhodd = 1/2$ and $\soc = 2$.
\epn

\bpf
The hexagonal grid is a 3-regular bipartite graph,
 thus both of its vertex classes are odd independent.
Moreover, the graph admits a perfect matching, hence independent sets cannot be denser than 1/2, also if the parity constraint of odd indepence is dropped.
\epf

\section{Concluding remarks and open problems}
\label{s:concl}

In this section we offer a list of open problems for future research.

\subsection{Grids}

\bpm
Determine $\sid$ and $\soc$ for the various types of $d$-dimensional grid
 graphs (Cartesian products of $d$ paths and/or cycles) for $d\geq 2$,
 or improve upon the bounds obtained in this research.
\epm

In  view of Theorem \ref{t:grid-d} we offer the following problem.

\bcj
The odd independence density of
 the $d$-dimensional grid $\gdinf = P_\infty \Boxx \cdots \Box
 P_\infty$ satisfies $\displaystyle \lim_{d\to\infty} \rhodd(\gdinf) = 1/3$.
\ecj

In the planar grid we constructed an odd-independent set of density 3/8, in which every vertex has 0 or 1 or 3 neighbors.
Perhaps this is essentially tight even if vertices with exactly two neighbors in the set are not excluded.

\bpm
Does there exist a sequence $\epsilon_n$ with $\lim_{n\to\infty} \epsilon_n = 0$,
 such that any $(\frac{3}{8}+\epsilon_n)\,n^2$ independent vertices of $P_n \Boxx P_n$
 contain all the four neighbors of some vertex?
\epm

If the answer is affirmative, then this would
 solve the odd-independence number problem of the infinite grid directly.

There does not exist an analogue of triangular and hexagonal grids for regular pentagons.
On the other hand, regular tilings of the plane with pentagons have been widely investigated, as surveyed nicely on the web page \cite{wiki:penta}.
For instance, it is known that there exist 15 types of convex pentagons which can tile the Euclidean plane in a periodic way.
There also exist monohedral pentagonal tilings with rotational symmetry.

\bpm
Analyze $\rhodd$ and $\soc$ of the graphs of the various pentagonal tilings of the plane.
\epm

Some of those tiligs are 3-regular, others also contain vertices of higher degree.
For example, in the ``Floret pentagonal tiling'' (an instance of types 1, 5 and 6, with a 6-tile primitive unit) the density of vertices of degree 6 is 1/9, and the density of degree-3 vertices all of whose neighbors have degree 3 is 2/3.
We can observe that its graph has $\rhodd\geq 1/3$, but we do not know how tight this estimate is.

\subsection{Chessboards}

Contrary to variants of the Rook graph,
 an extension of Proposition \ref{p:bishop} to
 non-square Bishop graphs seems to be more complicated.

\bpm
Determine $\sid$ and $\soc$ in the $p\times q$ analogue of the Bishop graph.
\epm

In view of Proposition \ref{p:knight} the following problem arises.

\bpm
Determine $\sid$ and $\soc$ of the Knight graph of size $n$, and more generally on the $p\times q$ board.
\epm

With regard to Propositions \ref{p:queen-10} and \ref{p:queen-inf} on the Queen graphs, we pose the following conjecture and open problem.

\bcj
The Queen graph of any finite size has $\sid=1$.
\ecj

\bpm
Decide whether the infinite Queen graph has $\sid=\infty$ or $\sid=1$
\epm

It is clear that any proper coloring of the infinite Queen graph
 (henc also a strong odd coloring)
 requires an infinite number of colors, as
 every line (horizontal, vertical, or diagonal)
 induces an infinite clique.
Yet it is not excluded---although very hard to believe---that
 there exists an infinite odd independent set.

For a detailed study of (standard) domination and independence
 in Bishop chessboard graphs in the plane and on surfaces
 we refer to the thesis \cite{thesis}.

\paragraph{Acknowledgement.}

This research was supported in part
 by the Slovenian Research Agency ARRS program P1-0383, projects J1-3002 and J1-4351, and the annual work program of Rudolfovo (R. \v S.) and
 by the ERC Advanced Grant ``ERMiD'' (Zs.~T.).

\newpage

\section{Appendix: Small grids}

\newcommand{\opti}[4]{\nin \underline{$k=#1$\,:}\quad $\max |S| = #2$\,;\quad $\max |S| /k^2 = #3$\,; \par $S$ = $\{$ #4 $\}$}

\newcommand{\optis}[4]{\nin \underline{$k=#1$\,:}\quad $\max |S| = #2$\,;\quad $\max |S| /k^2 \approx #3$\,; \par $S$ = $\{$ #4 $\}$}


In the first three subsections below, we list odd independent sets
 of sizes given in Table \ref{tab:grid}, hence verifying that the values
 of the table are lower bounds on $\sid$.
As already mentioned in the text, it has been checked by a computer code
 that no larger \sois s exist.
In the last subsection largest internally odd-independent sets are presented,
 with sizes given in Table \ref{tab:intern-odd}, that are
 relevant with respect to the infinite planar grid.

Throughout, the vertices of the grids are represented as ordered pairs
 $(i,j)$ of integers, $0\leq i,j\leq k-1$.
In particular, in internally odd-independent sets the odd neighborhood
 condition with respect to $S$ means a restriction on the vertices
 $(i,j)$ with $1\leq i,j\leq k-2$.

\subsection{Planar grids}

	\optis{3}{5}{0.555556}{(0, 0) (0, 2) (1, 1) (2, 0) (2, 2)}

\msk

\opti{4}{5}{0.3125}{(0, 1) (0, 3) (1, 2) (2, 1) (2, 3)}

\msk

\opti{5}{12}{0.48}{(0, 0) (0, 2) (0, 4) (1, 1) (1, 3) (2, 0) (2, 4) (3, 1) (3, 3) (4, 0) (4,~2) (4, 4)}

\msk

\optis{6}{12}{0.333333}{(1, 1) (1, 3) (1, 5) (2, 2) (2, 4) (3, 1) (3, 5) (4, 2) (4, 4) (5, 1) (5,~3) (5, 5)}

\msk

\optis{7}{20}{0.408163}{(0, 0) (0, 2) (0, 4) (0, 6) (1, 1) (1, 3) (1, 5) (2, 0) (2, 6) (3, 1) (3,~5) (4, 0) (4, 6) (5, 1) (5, 3) (5, 5) (6, 0) (6, 2) (6, 4) (6, 6)}

\msk

\opti{8}{21}{0.328125}{(0, 7) (1, 1) (1, 3) (1, 5) (2, 2) (2, 4) (3, 1) (3, 5) (3, 7) (4, 2) (4,~4) (4, 6) (5, 1) (5, 3) (5, 7) (6, 4) (6, 6) (7, 0) (7, 3) (7, 5) (7, 7)}

\msk

\optis{9}{29}{0.358025}{(0, 0) (0, 2) (0, 4) (0, 6) (0, 8) (1, 1) (1, 3) (1, 5) (1, 7) (2, 0) (2,~8) (3, 1) (3, 4) (3, 7) (4, 0) (4, 8) (5, 1) (5, 5) (5, 7) (6, 0) (6, 6) (6, 8) (7, 1) (7, 3) (7, 5) (8, 0) (8, 2) (8, 4) (8, 6)}

\msk

\opti{10}{33}{0.33}{(0, 0) (0, 2) (0, 4) (0, 7) (0, 9) (1, 1) (1, 3) (1, 8) (2, 0) (2, 4) (2,~7) (2, 9) (3, 1) (3, 3) (4, 0) (4, 2) (4, 4) (4, 6) (4, 8) (5, 5) (5, 7) (6, 4) (6, 8) (7, 0) (7, 2) (7, 5) (7, 7) (8, 1) (8, 4) (8, 6) (8, 8) (9, 0) (9, 2)}

\msk

\optis{11}{42}{0.347107}{(0, 0) (0, 2) (0, 6) (0, 8) (0, 10) (1, 1) (1, 7) (1, 9) (2, 0) (2, 2) (2,~4) (2, 6) (2, 10) (3, 3) (3, 5) (3, 7) (3, 9) (4, 2) (4, 8) (4, 10) (5, 3) (5, 7) (6, 0) (6, 2) (6, 8) (7, 1) (7, 3) (7, 5) (7, 7) (8, 0) (8, 4) (8, 6) (8, 8) (8, 10) (9, 1) (9, 3) (9, 9) (10, 0) (10, 2) (10, 4) (10, 8)}

\msk

\optis{12}{48}{0.333333}{(0, 0) (0, 2) (0, 4) (0, 7) (0, 9) (0, 11) (1, 1) (1, 3) (1, 8) (1, 10) (2, 0) (2, 4) (2, 7) (2, 11) (3, 1) (3, 3) (3, 8) (3, 10) (4, 0) (4, 2) (4, 4) (4, 7) (4, 9) (4, 11) (6, 6) (6, 8) (6, 10) (7, 0) (7, 2) (7, 4) (7, 7) (7, 9) (8, 1) (8, 3) (8, 6) (8, 10) (9, 0) (9, 4) (9, 7) (9, 9) (10, 1) (10, 3) (10, 6) (10, 8) (10, 10) (11, 0) (11, 2) (11, 4)}

\msk

\optis{13}{60}{0.35503}{(0, 0) (0, 2) (0, 4) (0, 8) (0, 10) (0, 12) (1, 1) (1, 3) (1, 6) (1, 9) (1,~11) (2, 0) (2, 4) (2, 8) (2, 12) (3, 1) (3, 3) (3, 9) (3, 11) (4, 0) (4, 2) (4,~4) (4, 6) (4, 8) (4, 10) (4, 12) (5, 5) (5, 7) (6, 1) (6, 4) (6, 8) (6, 11) (7, 5) (7,~7) (8, 0) (8, 2) (8, 4) (8, 6) (8, 8) (8, 10) (8, 12) (9, 1) (9, 3) (9, 9) (9, 11) (10,~0) (10, 4) (10, 8) (10, 12) (11, 1) (11, 3) (11, 6) (11, 9) (11, 11) (12, 0) (12, 2) (12, 4) (12, 8) (12, 10) (12, 12)}

\msk

\optis{14}{64}{0.326531}{(0, 6) (0, 8) (0, 10) (0, 13) (1, 0) (1, 2) (1, 4) (1, 7) (1, 9) (2, 1) (2,~3) (2, 6) (2, 10) (2, 12) (3, 0) (3, 4) (3, 7) (3, 9) (3, 11) (4, 1) (4, 3) (4,~6) (4, 8) (4, 12) (5, 0) (5, 2) (5, 4) (5, 9) (5, 11) (6, 8) (6, 10) (6, 12) (7, 1) (7,~3) (7, 5) (8, 2) (8, 4) (8, 9) (8, 11) (8, 13) (9, 1) (9, 5) (9, 7) (9, 10) (9,~12) (10,~2) (10, 4) (10, 6) (10, 9) (10, 13) (11, 1) (11, 3) (11, 7) (11, 10) (11, 12) (12, 4) (12, 6) (12, 9) (12, 11) (12, 13) (13, 0) (13, 3) (13, 5) (13, 7)}

\msk

\optis{15}{79}{0.351111}{(0, 0) (0, 2) (0, 4) (0, 8) (0, 10) (0, 12) (1, 1) (1, 3) (1, 6) (1, 9) (1,~11) (2, 0) (2, 4) (2, 8) (2, 12) (2, 14) (3, 1) (3, 3) (3, 9) (3, 11) (3, 13) (4,~0) (4, 2) (4, 4) (4, 6) (4, 8) (4, 10) (4, 14) (5, 5) (5, 7) (5, 11) (5, 13) (6,~1) (6,~4) (6, 8) (6, 10) (6, 12) (6, 14) (7, 5) (7, 7) (7, 9) (8, 0) (8, 2) (8,~4) (8,~6) (8, 10) (8, 13) (9, 1) (9, 3) (9, 7) (9, 9) (10, 0) (10, 4) (10, 6) (10,~8) (10,~10) (10, 12) (10, 14) (11, 1) (11, 3) (11, 5) (11, 11) (11, 13) (12, 0) (12,~2) (12,~6) (12, 10) (12, 14) (13, 3) (13, 5) (13, 8) (13, 11) (13, 13) (14, 2) (14, 4) (14,~6) (14, 10) (14, 12) (14, 14)}

\subsection{Cylindric grids}

Here we imagine the path $P_k$ horizontally and the cycle $C_k$ vertically.
That is, $(i,0)$ and $(i,k-1)$ are adjacent, while $(0,j)$ and $(k-1,j)$ are nonadjacent for all $0\leq i,j\leq k-1$.

\bsk

\optis{3}{2}{0.222222}{(0, 0) (2, 1)}

\msk

\opti{4}{5}{0.3125}{(0, 1) (0, 3) (1, 0) (1, 2) (3, 1)}

\msk

\opti{5}{6}{0.24}{(0, 0) (0, 3) (1, 4) (2, 0) (2, 3) (4, 1)}

\msk

\optis{6}{12}{0.333333}{(0, 0) (0, 2) (0, 4) (1, 1) (1, 3) (1, 5) (4, 0) (4, 2) (4, 4) (5, 1) (5,~3) (5, 5)}

\msk

\optis{7}{14}{0.285714}{(0, 0) (0, 3) (0, 5) (1, 4) (1, 6) (2, 0) (2, 3) (3, 4) (3, 6) (4, 0) (4,~3) (4, 5) (6, 1) (6, 4)}

\msk

\opti{8}{24}{0.375}{(0, 1) (0, 3) (0, 5) (0, 7) (1, 0) (1, 2) (1, 4) (1, 6) (3, 1) (3, 3) (3,~5) (3, 7) (4, 0) (4, 2) (4, 4) (4, 6) (6, 1) (6, 3) (6, 5) (6, 7) (7, 0) (7, 2) (7, 4) (7, 6)}

\msk

\optis{9}{25}{0.308642}{(0, 1) (0, 3) (0, 6) (0, 8) (1, 0) (1, 2) (1, 7) (2, 3) (2, 6) (3, 2) (3,~7) (4, 0) (4, 3) (4, 6) (5, 2) (5, 7) (6, 3) (6, 6) (7, 0) (7, 2) (7, 7) (8, 1) (8, 3) (8, 6) (8, 8)}

\msk

\opti{10}{34}{0.34}{(0, 0) (0, 2) (0, 4) (0, 6) (0, 8) (1, 1) (1, 3) (1, 5) (1, 7) (1, 9) (3,~0) (3, 2) (3, 4) (3, 6) (3, 8) (4, 1) (4, 3) (4, 5) (4, 7) (4, 9) (6, 0) (6, 2) (6, 6) (6, 8) (7, 1) (7, 7) (7, 9) (8, 2) (8, 4) (8, 6) (9, 1) (9, 3) (9, 5) (9, 7)}

\msk

\optis{11}{37}{0.305785}{(0, 1) (0, 4) (0, 6) (0, 8) (1, 5) (1, 7) (1, 10) (2, 2) (2, 4) (2, 8) (3,~3) (3, 5) (3, 7) (4, 0) (4, 2) (4, 6) (4, 8) (5, 1) (5, 3) (5, 5) (6, 0) (6, 4) (6,p~6) (6, 9) (7, 1) (7, 3) (7, 10) (8, 2) (8, 4) (8, 9) (9, 1) (9, 7) (9, 10) (10,~0) (10, 2) (10, 5) (10, 9)}

\msk

\optis{12}{48}{0.333333}{(0, 0) (0, 2) (0, 4) (0, 6) (0, 8) (0, 10) (1, 1) (1, 3) (1, 5) (1, 7) (1,~9) (1, 11) (3, 0) (3, 2) (3, 4) (3, 6) (3, 8) (3, 10) (4, 1) (4, 3) (4, 5) (4, 7) (4, 9) (4, 11) (6, 6) (6, 8) (6, 10) (7, 0) (7, 2) (7, 4) (7, 7) (7, 9) (8, 1) (8, 3) (8, 6) (8, 10) (9, 0) (9, 4) (9, 7) (9, 9) (10, 1) (10, 3) (10, 6) (10, 8) (10, 10) (11, 0) (11, 2) (11, 4)}

\msk

\optis{13}{52}{0.307692}{(0, 0) (0, 2) (0, 4) (0, 6) (0, 8) (0, 10) (1, 1) (1, 3) (1, 5) (1, 7) (1,~9) (2, 0) (2, 10) (3, 1) (3, 9) (4, 0) (4, 3) (4, 5) (4, 7) (4, 10) (5, 1) (5, 4) (5, 6) (5, 9) (6, 0) (6, 3) (6, 7) (6, 10) (7, 1) (7, 4) (7, 6) (7, 9) (8, 0) (8, 3) (8, 5) (8, 7) (8, 10) (9, 1) (9, 9) (10, 0) (10, 10) (11, 1) (11, 3) (11, 5) (11, 7) (11, 9) (12, 0) (12, 2) (12, 4) (12, 6) (12, 8) (12, 10)}

\msk

\optis{14}{70}{0.357143}{(0, 1) (0, 3) (0, 5) (0, 7) (0, 9) (0, 11) (0, 13) (1, 0) (1, 2) (1, 4) (1,~6) (1, 8) (1, 10) (1, 12) (3, 1) (3, 3) (3, 5) (3, 7) (3, 9) (3, 11) (3, 13) (4,~0) (4, 2) (4, 4) (4, 6) (4, 8) (4, 10) (4, 12) (6, 1) (6, 3) (6, 5) (6, 7) (6,~9) (6, 11) (6, 13) (7, 0) (7, 2) (7, 4) (7, 6) (7, 8) (7, 10) (7, 12) (9, 1) (9, 3) (9,~5) (9,~7) (9, 9) (9, 11) (9, 13) (10, 0) (10, 2) (10, 4) (10, 6) (10, 8) (10,~10) (10,~12) (12, 1) (12, 3) (12, 5) (12, 7) (12, 9) (12, 11) (12, 13) (13, 0) (13, 2) (13, 4) (13, 6) (13, 8) (13, 10) (13, 12)}

\subsection{Toroidal grids}

\optis{3}{1}{0.111111}{(1, 2)}

\msk

\opti{4}{6}{0.375}{(0, 1) (0, 3) (1, 2) (2, 1) (2, 3) (3, 0)}

\msk

\opti{5}{5}{0.2}{(2, 1) (2, 3) (3, 2) (4, 1) (4, 3)}

\msk

\optis{6}{12}{0.333333}{(0, 2) (0, 5) (1, 0) (1, 3) (2, 2) (2, 5) (3, 0) (3, 3) (4, 2) (4, 5) (5,~0) (5, 3)}

\msk

\optis{7}{12}{0.244898}{(0, 0) (0, 3) (0, 5) (3, 0) (3, 3) (3, 5) (4, 4) (4, 6) (5, 0) (5, 3) (6,~4) (6, 6)}

\msk

\opti{8}{24}{0.375}{(0, 0) (0, 2) (0, 4) (0, 6) (1, 3) (1, 7) (2, 0) (2, 2) (2, 4) (2, 6) (3,~1) (3, 5) (4, 0) (4, 2) (4, 4) (4, 6) (5, 3) (5, 7) (6, 0) (6, 2) (6, 4) (6, 6) (7, 1) (7, 5)}

\msk

\optis{9}{21}{0.259259}{(0, 2) (0, 4) (1, 6) (2, 1) (2, 8) (3, 0) (3, 3) (3, 5) (4, 1) (4, 4) (4,~8) (5, 3) (5, 5) (6, 0) (6, 7) (7, 2) (7, 4) (7, 8) (8, 0) (8, 3) (8, 7)}

\msk

\opti{10}{30}{0.3}{(0, 1) (0, 3) (0, 5) (0, 7) (0, 9) (1, 0) (1, 8) (2, 1) (2, 4) (2, 7) (3,~0) (3, 8) (4, 1) (4, 3) (4, 5) (4, 7) (4, 9) (5, 2) (5, 4) (5, 6) (6, 1) (6, 7) (7, 2) (7, 6) (7, 9) (8, 1) (8, 7) (9, 2) (9, 4) (9, 6)}

\msk

\optis{11}{34}{0.280992}{(0, 0) (0, 2) (0, 4) (1, 1) (1, 5) (1, 8) (1, 10) (2, 2) (2, 4) (2, 9) (3,~1) (3, 3) (3, 5) (3, 8) (3, 10) (4, 0) (5, 1) (5, 10) (6, 3) (6, 5) (6, 8) (7,~4) (8, 1) (8, 3) (8, 5) (8, 7) (8, 10) (9, 0) (9, 2) (9, 6) (10, 3) (10, 5) (10, 7) (10,~10)}

\msk

\opti{12}{54}{0.375}{(0, 2) (0, 4) (0, 6) (0, 8) (0, 10) (1, 1) (1, 5) (1, 9) (1, 11) (2, 0) (2,~2) (2, 4) (2, 6) (2, 8) (3, 3) (3, 7) (3, 9) (3, 11) (4, 0) (4, 2) (4, 4) (4, 6) (4, 10) (5, 1) (5, 5) (5, 7) (5, 9) (6, 0) (6, 2) (6, 4) (6, 8) (6, 10) (7, 3) (7, 5) (7, 7) (7, 11) (8, 0) (8, 2) (8, 6) (8, 8) (8, 10) (9, 1) (9, 3) (9, 5) (9, 9) (10,~0) (10, 4) (10, 6) (10, 8) (10, 10) (11, 1) (11, 3) (11, 7) (11, 11)}

\msk

\optis{13}{49}{0.289941}{(0, 2) (0, 4) (0, 8) (0, 10) (0, 12) (1, 1) (1, 3) (1, 5) (2, 7) (2, 9) (2,~11) (3, 0) (3, 2) (3, 4) (3, 8) (3, 10) (4, 1) (4, 3) (4, 7) (4, 11) (5, 0) (5,~4) (5, 8) (5, 10) (6, 1) (6, 3) (6, 7) (6, 9) (6, 11) (7, 0) (7, 2) (7, 4) (8, 6) (9, 8) (9, 10) (9, 12) (10, 1) (10, 3) (10, 5) (10, 9) (10, 11) (11, 2) (11, 4) (11, 8) (11, 12) (12, 1) (12, 5) (12, 9) (12, 11)}

\msk

\optis{14}{62}{0.316327}{(0, 1) (0, 3) (0, 5) (0, 8) (0, 10) (1, 2) (1, 4) (1, 7) (1, 11) (1, 13) (2, 1) (2, 5) (2, 8) (2, 10) (2, 12) (3, 2) (3, 4) (3, 7) (3, 9) (3, 13) (4, 1) (4,~3) (4, 5) (4, 10) (4, 12) (5, 9) (5, 11) (5, 13) (6, 2) (6, 4) (6, 6) (7, 3) (7, 5) (7,~8) (7, 10) (7, 12) (8, 0) (8, 2) (8, 6) (8, 9) (8, 11) (9, 1) (9, 3) (9, 5) (9,~8) (9, 12) (10, 0) (10, 4) (10, 6) (10, 9) (10, 11) (11, 1) (11, 3) (11, 8) (11, 10) (11, 12) (12, 0) (12, 2) (12, 4) (13, 7) (13, 9) (13, 11)}

\subsection{Largest internally odd-independent sets}

\opti{5}{12}{0.48}{(0, 0) (0, 2) (0, 4) (1, 1) (1, 3) (2, 0) (2, 4) (3, 1) (3, 3) (4, 0) (4,~2) (4, 4)}

\msk

\optis{6}{15}{0.416667}{(0, 1) (0, 3) (0, 5) (1, 0) (1, 4) (2, 1) (2, 3) (2, 5) (3, 0) (3, 2) (4,~3) (4, 5) (5, 0) (5, 2) (5, 4)}

\msk

\optis{7}{21}{0.428571}{(0, 0) (0, 2) (0, 4) (0, 6) (1, 1) (1, 5) (2, 0) (2, 2) (2, 4) (2, 6) (3,~3) (4, 0) (4, 2) (4, 4) (4, 6) (5, 1) (5, 5) (6, 0) (6, 2) (6, 4) (6, 6)}

\msk

\opti{8}{26}{0.40625}{(0, 0) (0, 2) (0, 4) (0, 6) (1, 3) (1, 7) (2, 0) (2, 2) (2, 4) (2, 6) (3,~1) (3, 5) (3, 7) (4, 0) (4, 2) (4, 4) (5, 3) (5, 5) (5, 7) (6, 0) (6, 2) (6, 6) (7, 1) (7, 3) (7, 5) (7, 7)}

\msk

\optis{9}{34}{0.419753}{(0, 0) (0, 2) (0, 4) (0, 6) (0, 8) (1, 1) (1, 5) (1, 7) (2, 0) (2, 2) (2,~4) (2, 8) (3, 3) (3, 5) (3, 7) (4, 0) (4, 2) (4, 6) (4, 8) (5, 1) (5, 3) (5, 5) (6, 0) (6, 4) (6, 6) (6, 8) (7, 1) (7, 3) (7, 7) (8, 0) (8, 2) (8, 4) (8, 6) (8, 8)}

\msk

\opti{10}{40}{0.4}{(0, 1) (0, 3) (0, 5) (0, 7) (0, 9) (1, 0) (1, 4) (1, 8) (2, 1) (2, 3) (2,~5) (2, 7) (2, 9) (3, 0) (3, 2) (3, 6) (4, 3) (4, 5) (4, 7) (4, 9) (5, 0) (5, 2) (5, 4) (5, 8) (6, 1) (6, 5) (6, 7) (6, 9) (7, 0) (7, 2) (7, 4) (7, 6) (8, 3) (8, 7) (8, 9) (9, 0) (9, 2) (9, 4) (9, 6) (9, 8)}

\msk

\optis{11}{49}{0.404959}{(0, 0) (0, 2) (0, 4) (0, 6) (0, 8) (0, 10) (1, 3) (1, 7) (1, 9) (2, 0) (2,~2) (2, 4) (2, 6) (2, 10) (3, 1) (3, 5) (3, 7) (3, 9) (4, 0) (4, 2) (4, 4) (4,~8) (4,~10) (5, 3) (5, 5) (5, 7) (6, 0) (6, 2) (6, 6) (6, 8) (6, 10) (7, 1) (7, 3) (7,~5) (7, 9) (8, 0) (8, 4) (8, 6) (8, 8) (8, 10) (9, 1) (9, 3) (9, 7) (10, 0) (10, 2) (10,~4) (10, 6) (10, 8) (10, 10)}

\msk

\optis{12}{57}{0.395833}{(0, 0) (0, 2) (0, 4) (0, 6) (0, 8) (0, 10) (1, 1) (1, 5) (1, 9) (1, 11) (2,~0) (2, 2) (2, 4) (2, 6) (2, 8) (3, 3) (3, 7) (3, 9) (3, 11) (4, 0) (4, 2) (4,~4) (4,~6) (4, 10) (5, 1) (5, 5) (5, 7) (5, 9) (5, 11) (6, 0) (6, 2) (6, 4) (6, 8) (7,~3) (7, 5) (7, 7) (7, 9) (7, 11) (8, 0) (8, 2) (8, 6) (8, 10) (9, 1) (9, 3) (9, 5) (9,~7) (9, 9) (9, 11) (10, 0) (10, 4) (10, 8) (11, 1) (11, 3) (11, 5) (11, 7) (11, 9) (11,~11)}

\msk

\optis{13}{68}{0.402367}{(0, 0) (0, 2) (0, 4) (0, 6) (0, 8) (0, 10) (0, 12) (1, 1) (1, 5) (1, 7) (1,~11) (2, 0) (2, 2) (2, 4) (2, 8) (2, 10) (2, 12) (3, 3) (3, 5) (3, 7) (3, 9) (4,~0) (4, 2) (4, 6) (4, 10) (4, 12) (5, 1) (5, 3) (5, 5) (5, 7) (5, 9) (5, 11) (6, 0) (6,~4) (6, 8) (6, 12) (7, 1) (7, 3) (7, 5) (7, 7) (7, 9) (7, 11) (8, 0) (8, 2) (8, 6) (8,~10) (8, 12) (9, 3) (9, 5) (9, 7) (9, 9) (10, 0) (10, 2) (10, 4) (10, 8) (10, 10) (10,~12) (11, 1) (11, 5) (11, 7) (11, 11) (12, 0) (12, 2) (12, 4) (12, 6) (12, 8) (12, 10) (12, 12)}

\msk

\optis{14}{77}{0.392857}{(0, 0) (0, 2) (0, 4) (0, 6) (0, 8) (0, 10) (0, 12) (1, 1) (1, 5) (1, 9) (1,~13) (2, 0) (2, 2) (2, 4) (2, 6) (2, 8) (2, 10) (2, 12) (3, 3) (3, 7) (3, 11) (3,~13) (4, 0) (4, 2) (4, 4) (4, 6) (4, 8) (4, 10) (5, 1) (5, 5) (5, 9) (5, 11) (5,~13) (6, 0) (6, 2) (6, 4) (6, 6) (6, 8) (6, 12) (7, 3) (7, 7) (7, 9) (7, 11) (7,~13) (8,~0) (8, 2) (8, 4) (8, 6) (8, 10) (9, 1) (9, 5) (9, 7) (9, 9) (9, 11) (9, 13) (10, 0) (10,~2) (10, 4) (10, 8) (10, 12) (11, 3) (11, 5) (11, 7) (11, 9) (11, 11) (11, 13) (12, 0) (12, 2) (12, 6) (12, 10) (13, 1) (13, 3) (13, 5) (13, 7) (13, 9) (13, 11) (13, 13)}

\msk

\optis{15}{89}{0.395556}{(0, 0) (0, 2) (0, 4) (0, 6) (0, 8) (0, 10) (0, 12) (0, 14) (1, 1) (1, 5) (1,~9) (1, 13) (2, 0) (2, 2) (2, 4) (2, 6) (2, 8) (2, 10) (2, 12) (2, 14) (3, 3) (3,~7) (3, 11) (4, 0) (4, 2) (4, 4) (4, 6) (4, 8) (4, 10) (4, 12) (4, 14) (5, 1) (5,~5) (5,~9) (5, 13) (6, 0) (6, 2) (6, 4) (6, 6) (6, 8) (6, 10) (6, 12) (6, 14) (7,~3) (7,~7) (7,~11) (8, 0) (8, 2) (8, 4) (8, 6) (8, 8) (8, 10) (8, 12) (8, 14) (9,~1) (9,~5) (9,~9) (9, 13) (10, 0) (10, 2) (10, 4) (10, 6) (10, 8) (10, 10) (10,~12) (10,~14) (11, 3) (11, 7) (11, 11) (12, 0) (12, 2) (12, 4) (12, 6) (12, 8) (12, 10) (12, 12) (12, 14) (13, 1) (13, 5) (13, 9) (13, 13) (14, 0) (14, 2) (14, 4) (14, 6) (14, 8) (14, 10) (14, 12) (14, 14)}

\msk

\opti{16}{100}{0.390625}{(0, 0) (0, 2) (0, 4) (0, 6) (0, 8) (0, 10) (0, 12) (0, 14) (1, 1) (1,~5) (1,~9) (1, 13) (1, 15) (2, 0) (2, 2) (2, 4) (2, 6) (2, 8) (2, 10) (2, 12) (3, 3) (3,~7) (3, 11) (3, 13) (3, 15) (4, 0) (4, 2) (4, 4) (4, 6) (4, 8) (4, 10) (4, 14) (5,~1) (5,~5) (5, 9) (5, 11) (5, 13) (5, 15) (6, 0) (6, 2) (6, 4) (6, 6) (6, 8) (6,~12) (7,~3) (7, 7) (7, 9) (7, 11) (7, 13) (7, 15) (8, 0) (8, 2) (8, 4) (8, 6) (8, 10) (8,~14) (9, 1) (9, 5) (9, 7) (9, 9) (9, 11) (9, 13) (9, 15) (10, 0) (10, 2) (10, 4) (10, 8) (10, 12) (11, 3) (11, 5) (11, 7) (11, 9) (11, 11) (11, 13) (11, 15) (12, 0) (12,~2) (12, 6) (12, 10) (12, 14) (13, 1) (13, 3) (13, 5) (13, 7) (13, 9) (13, 11) (13, 13) (13, 15) (14, 0) (14, 4) (14, 8) (14, 12) (15, 1) (15, 3) (15, 5) (15, 7) (15, 9) (15, 11) (15, 13) (15, 15)}

\msk

\optis{17}{114}{0.394464}{(0, 0) (0, 2) (0, 4) (0, 6) (0, 8) (0, 10) (0, 12) (0, 14) (0, 16) (1, 1) (1, 3) (1, 7) (1, 11) (1, 15) (2, 0) (2, 4) (2, 6) (2, 8) (2, 10) (2, 12) (2,~14) (2,~16) (3, 1) (3, 3) (3, 5) (3, 9) (3, 13) (4, 0) (4, 2) (4, 6) (4, 8) (4, 10) (4,~12) (4, 14) (4, 16) (5, 3) (5, 5) (5, 7) (5, 11) (5, 15) (6, 0) (6, 2) (6, 4) (6,~8) (6,~10) (6, 12) (6, 14) (6, 16) (7, 1) (7, 5) (7, 7) (7, 9) (7, 13) (8, 0) (8,~2) (8,~4) (8, 6) (8, 10) (8, 12) (8, 14) (8, 16) (9, 3) (9, 7) (9, 9) (9, 11) (9,~15) (10, 0) (10, 2) (10, 4) (10, 6) (10, 8) (10, 12) (10, 14) (10, 16) (11,~1) (11, 5) (11, 9) (11, 11) (11, 13) (12, 0) (12, 2) (12, 4) (12, 6) (12, 8) (12,~10) (12,~14) (12, 16) (13, 3) (13, 7) (13, 11) (13, 13) (13, 15) (14, 0) (14, 2) (14,~4) (14,~6) (14, 8) (14, 10) (14, 12) (14, 16) (15, 1) (15, 5) (15, 9) (15, 13) (15,~15) (16,~0) (16, 2) (16, 4) (16, 6) (16, 8) (16, 10) (16, 12) (16, 14) (16, 16)}

\msk

\optis{18}{125}{0.388889}{(0, 0) (0, 2) (0, 4) (0, 6) (0, 8) (0, 10) (0, 12) (0, 14) (0, 16) (1,~3) (1, 7) (1, 11) (1, 15) (1, 17) (2, 0) (2, 2) (2, 4) (2, 6) (2, 8) (2, 10) (2, 12) (2,~14) (3, 1) (3, 5) (3, 9) (3, 13) (3, 15) (3, 17) (4, 0) (4, 2) (4, 4) (4, 6) (4,~8) (4, 10) (4, 12) (4, 16) (5, 3) (5, 7) (5, 11) (5, 13) (5, 15) (5, 17) (6,~0) (6, 2) (6,~4) (6, 6) (6, 8) (6, 10) (6, 14) (7, 1) (7, 5) (7, 9) (7, 11) (7, 13) (7,~15) (7,~17) (8, 0) (8, 2) (8, 4) (8, 6) (8, 8) (8, 12) (8, 16) (9, 3) (9, 7) (9,~9) (9,~11) (9, 13) (9, 15) (9, 17) (10, 0) (10, 2) (10, 4) (10, 6) (10, 10) (10,~14) (11, 1) (11, 5) (11, 7) (11, 9) (11, 11) (11, 13) (11, 15) (11, 17) (12,~0) (12,~2) (12, 4) (12, 8) (12, 12) (12, 16) (13, 3) (13, 5) (13, 7) (13,~9) (13,~11) (13,~13) (13,~15) (13, 17) (14, 0) (14, 2) (14, 6) (14, 10) (14, 14) (15,~1) (15,~3) (15,~5) (15,~7) (15, 9) (15, 11) (15, 13) (15, 15) (15, 17) (16,~0) (16, 4) (16, 8) (16,~12) (16,~16) (17, 1) (17, 3) (17, 5) (17, 7) (17, 9) (17, 11) (17, 13) (17,~15) (17,~17)}

\msk

\optis{19}{141}{0.390582}{(0, 0) (0, 2) (0, 4) (0, 6) (0, 8) (0, 10) (0, 12) (0, 14) (0, 16) (0,~18) (1, 3) (1, 7) (1, 9) (1, 13) (1, 17) (2, 0) (2, 2) (2, 4) (2, 6) (2, 10) (2, 12) (2,~14) (2, 16) (2, 18) (3, 1) (3, 5) (3, 7) (3, 9) (3, 11) (3, 15) (4,~0) (4,~2) (4,~4) (4, 8) (4, 12) (4, 14) (4, 16) (4, 18) (5, 3) (5, 5) (5, 7) (5, 9) (5, 11) (5,~13) (5, 17) (6,~0) (6, 2) (6, 6) (6, 10) (6, 14) (6, 16) (6, 18) (7,~1) (7,~3) (7,~5) (7,~7) (7,~9) (7, 11) (7, 13) (7, 15) (8, 0) (8, 4) (8, 8) (8, 12) (8, 16) (8,~18) (9, 1) (9,~3) (9,~5) (9, 7) (9, 9) (9, 11) (9, 13) (9, 15) (9,~17) (10,~0) (10,~2) (10,~6) (10,~10) (10, 14) (10, 18) (11, 3) (11, 5) (11, 7) (11,~9) (11,~11) (11,~13) (11,~15) (11,~17) (12, 0) (12, 2) (12, 4) (12, 8) (12, 12) (12,~16) (12,~18) (13,~1) (13, 5) (13, 7) (13, 9) (13, 11) (13, 13) (13, 15) (14,~0) (14,~2) (14, 4) (14,~6) (14, 10) (14, 14) (14, 16) (14, 18) (15, 3) (15,~7) (15, 9) (15,~11) (15, 13) (15,~17) (16, 0) (16, 2) (16, 4) (16, 6) (16, 8) (16, 12) (16, 14) (16,~16) (16, 18) (17, 1) (17, 5) (17, 9) (17, 11) (17, 15) (18, 0) (18, 2) (18, 4) (18, 6) (18, 8) (18, 10) (18, 12) (18, 14) (18, 16) (18, 18)}

\msk

\opti{20}{155}{0.3875}{(0, 0) (0, 2) (0, 4) (0, 6) (0, 8) (0, 10) (0, 12) (0, 14) (0, 16) (0,~18) (1, 1) (1, 5) (1, 9) (1, 13) (1, 17) (1, 19) (2, 0) (2, 2) (2, 4) (2, 6) (2,~8) (2,~10) (2, 12) (2, 14) (2, 16) (3, 3) (3, 7) (3, 11) (3, 15) (3, 17) (3,~19) (4,~0) (4, 2) (4,~4) (4, 6) (4, 8) (4, 10) (4, 12) (4, 14) (4, 18) (5, 1) (5, 5) (5, 9) (5,~13) (5,~15) (5, 17) (5, 19) (6, 0) (6, 2) (6, 4) (6, 6) (6, 8) (6, 10) (6,~12) (6, 16) (7,~3) (7,~7) (7, 11) (7, 13) (7, 15) (7, 17) (7, 19) (8, 0) (8,~2) (8,~4) (8, 6) (8,~8) (8,~10) (8, 14) (8, 18) (9, 1) (9, 5) (9, 9) (9, 11) (9, 13) (9, 15) (9, 17) (9, 19) (10, 0) (10, 2) (10, 4) (10, 6) (10, 8) (10, 12) (10,~16) (11,~3) (11,~7) (11,~9) (11, 11) (11, 13) (11, 15) (11, 17) (11, 19) (12, 0) (12,~2) (12,~4) (12,~6) (12,~10) (12, 14) (12, 18) (13, 1) (13, 5) (13,~7) (13,~9) (13,~11) (13,~13) (13,~15) (13, 17) (13, 19) (14, 0) (14, 2) (14, 4) (14, 8) (14,~12) (14,~16) (15,~3) (15, 5) (15, 7) (15, 9) (15, 11) (15, 13) (15, 15) (15,~17) (15,~19) (16,~0) (16,~2) (16, 6) (16, 10) (16, 14) (16, 18) (17, 1) (17,~3) (17,~5) (17, 7) (17, 9) (17,~11) (17, 13) (17, 15) (17, 17) (17, 19) (18, 0) (18,~4) (18,~8) (18, 12) (18,~16) (19,~1) (19, 3) (19, 5) (19, 7) (19, 9) (19, 11) (19,~13) (19,~15) (19, 17)\break (19,~19)}

\msk

\opti{21}{172}{0.390023}{%
(0, 0) (0, 2) (0, 4) (0, 6) (0, 8) (0, 10) (0, 12) (0, 14) (0, 16) (0,~18) (0, 20) (1, 1) (1, 5) (1, 9) (1, 13) (1, 15) (1, 19) (2, 0) (2, 2) (2,~4) (2, 6) (2,~8) (2, 10) (2, 12) (2, 16) (2, 18) (2, 20) (3, 3) (3, 7) (3, 11) (3,~13) (3,~15) (3,~17) (4, 0) (4, 2) (4, 4) (4, 6) (4, 8) (4, 10) (4, 14) (4, 18) (4,~20) (5, 1) (5,~5) (5,~9) (5, 11) (5, 13) (5, 15) (5, 17) (5, 19) (6, 0) (6, 2) (6,~4) (6, 6) (6,~8) (6,~12) (6,~16) (6, 20) (7, 3) (7, 7) (7, 9) (7, 11) (7, 13) (7,~15) (7,~17) (7,~19) (8,~0) (8,~2) (8, 4) (8, 6) (8, 10) (8, 14) (8, 18) (8, 20) (9,~1) (9,~5) (9,~7) (9,~9) (9,~11) (9, 13) (9, 15) (9, 17) (10, 0) (10, 2) (10, 4) (10,~8) (10,~12) (10,~16) (10,~18) (10, 20) (11, 3) (11, 5) (11, 7) (11, 9) (11,~11) (11,~13) (11,~15) (11,~19) (12,~0) (12, 2) (12, 6) (12, 10) (12, 14) (12, 16) (12,~18) (12,~20) (13,~1) (13,~3) (13,~5) (13, 7) (13, 9) (13, 11) (13, 13) (13,~17) (14,~0) (14,~4) (14, 8) (14, 12) (14,~14) (14, 16) (14, 18) (14, 20) (15, 1) (15,~3) (15,~5) (15,~7) (15,~9) (15,~11) (15,~15) (15, 19) (16, 0) (16, 2) (16, 6) (16,~10) (16,~12) (16,~14) (16,~16) (16, 18) (16,~20) (17, 3) (17, 5) (17, 7) (17, 9) (17,~13) (17,~17) (18, 0) (18, 2) (18,~4) (18, 8) (18, 10) (18, 12) (18, 14) (18,~16) (18, 18) (18, 20) (19, 1) (19, 5) (19,~7) (19, 11) (19, 15) (19, 19) (20, 0) (20, 2) (20, 4) (20, 6) (20, 8) (20,~10) (20, 12) (20, 14) (20, 16) (20, 18) (20, 20)}

\msk

\opti{22}{187}{0.386364}{%
(0, 1) (0, 3) (0, 5) (0, 7) (0, 9) (0, 11) (0, 13) (0, 15) (0, 17) (0, 19) (0, 21) (1, 0) (1, 4) (1, 8) (1, 12) (1, 16) (1, 20) (2, 1) (2, 3) (2, 5) (2, 7) (2, 9) (2, 11) (2, 13) (2, 15) (2, 17) (2, 19) (2, 21) (3, 0) (3, 2) (3, 6) (3, 10) (3, 14) (3, 18) (4, 3) (4, 5) (4, 7) (4, 9) (4, 11) (4, 13) (4, 15) (4,~17) (4, 19) (4, 21) (5, 0) (5, 2) (5, 4) (5, 8) (5, 12) (5, 16) (5, 20) (6, 1) (6,~5) (6, 7) (6, 9) (6, 11) (6, 13) (6, 15) (6, 17) (6, 19) (6, 21) (7, 0) (7, 2) (7,~4) (7, 6) (7, 10) (7, 14) (7, 18) (8, 3) (8, 7) (8, 9) (8, 11) (8, 13) (8, 15) (8,~17) (8, 19) (8, 21) (9, 0) (9, 2) (9, 4) (9, 6) (9, 8) (9, 12) (9, 16) (9, 20) (10, 1) (10, 5) (10, 9) (10,~11) (10, 13) (10, 15) (10, 17) (10, 19) (10, 21) (11,~0) (11, 2) (11, 4) (11,~6) (11,~8) (11, 10) (11, 14) (11, 18) (12, 3) (12,~7) (12, 11) (12,~13) (12, 15) (12,~17) (12,~19) (12, 21) (13, 0) (13, 2) (13, 4) (13,~6) (13,~8) (13, 10) (13,~12) (13,~16) (13, 20) (14, 1) (14, 5) (14, 9) (14, 13) (14,~15) (14,~17) (14,~19) (14,~21) (15,~0) (15, 2) (15, 4) (15, 6) (15, 8) (15,~10) (15,~12) (15,~14) (15,~18) (16,~3) (16, 7) (16, 11) (16, 15) (16, 17) (16,~19) (16,~21) (17,~0) (17, 2) (17,~4) (17,~6) (17,~8) (17, 10) (17,~12) (17,~14) (17,~16) (17,~20) (18,~1) (18, 5) (18,~9) (18,~13) (18,~17) (18, 19) (18,~21) (19,~0) (19,~2) (19,~4) (19,~6) (19,~8) (19,~10) (19,~12) (19,~14) (19,~16) (19,~18) (20,~3) (20,~7) (20,~11) (20,~15) (20,~19) (20,~21) (21,~0) (21,~2) (21,~4) (21,~6) (21,~8) (21,~10) (21,~12) (21,~14) (21,~16) (21,~18) (21,~20)}

\msk

\opti{23}{204}{0.3875236294896006}{(0, 0) (0, 2) (0, 4) (0, 6) (0, 8) (0, 10) (0, 12) (0,~14) (0,~16) (0,~18) (0,~20) (0,~22) (1, 1) (1, 5) (1, 9) (1, 11) (1, 15) (1, 19) (2, 0) (2,~2) (2,~4) (2,~6) (2,~8) (2,~12) (2, 14) (2, 16) (2, 18) (2, 20) (2, 22) (3, 3) (3,~7) (3,~9) (3,~11) (3,~13) (3,~17) (3, 21) (4, 0) (4, 2) (4, 4) (4, 6) (4, 10) (4,~14) (4,~16) (4,~18) (4,~20) (4,~22) (5, 1) (5, 5) (5, 7) (5, 9) (5, 11) (5, 13) (5,~15) (5,~19) (6,~0) (6,~2) (6,~4) (6, 8) (6, 12) (6, 16) (6, 18) (6, 20) (6, 22) (7,~3) (7,~5) (7,~7) (7,~9) (7,~11) (7, 13) (7, 15) (7, 17) (7, 21) (8, 0) (8, 2) (8,~6) (8,~10) (8,~14) (8,~18) (8,~20) (8, 22) (9, 1) (9, 3) (9, 5) (9, 7) (9, 9) (9,~11) (9,~13) (9,~15) (9,~17) (9,~19) (10, 0) (10, 4) (10, 8) (10, 12) (10, 16) (10,~20) (10,~22) (11,~1) (11,~3) (11,~5) (11, 7) (11, 9) (11, 11) (11, 13) (11,~15) (11,~17) (11,~19) (11,~21) (12,~0) (12, 2) (12, 6) (12, 10) (12, 14) (12, 18) (12,~22) (13,~3) (13,~5) (13,~7) (13,~9) (13, 11) (13, 13) (13, 15) (13, 17) (13,~19) (13,~21) (14,~0) (14,~2) (14,~4) (14,~8) (14,~12) (14,~16) (14,~20) (14,~22) (15,~1) (15,~5) (15,~7) (15,~9) (15,~11) (15,~13) (15,~15) (15,~17) (15,~19) (16,~0) (16,~2) (16,~4) (16,~6) (16,~10) (16,~14) (16,~18) (16, 20) (16,~22) (17,~3) (17,~7) (17,~9) (17,~11) (17,~13) (17,~15) (17,~17) (17,~21) (18,~0) (18,~2) (18,~4) (18,~6) (18,~8) (18,~12) (18,~16) (18, 18) (18,~20) (18,~22) (19,~1) (19,~5) (19,~9) (19,~11) (19,~13) (19,~15) (19,~19) (20, 0) (20, 2) (20, 4) (20, 6) (20,~8) (20,~10) (20,~14) (20,~16) (20,~18) (20, 20) (20, 22) (21, 3) (21, 7) (21,~11) (21,~13) (21,~17) (21,~21) (22,~0) (22, 2) (22, 4) (22, 6) (22, 8) (22, 10) (22,~12) (22,~14) (22,~16) (22,~18) (22,~20) (22, 22)}

\msk

\opti{24}{222}{0.385417}{%
(0, 0) (0, 2) (0, 4) (0, 6) (0, 8) (0, 10) (0, 12) (0, 14) (0, 16) (0,~18) (0,~20) (0,~22) (1,~3) (1,~7) (1, 11) (1, 15) (1, 19) (1, 23) (2,~0) (2,~2) (2,~4) (2,~6) (2,~8) (2, 10) (2, 12) (2, 14) (2, 16) (2, 18) (2,~20) (2,~22) (3,~1) (3,~5) (3,~9) (3,~13) (3, 17) (3, 21) (3, 23) (4, 0) (4, 2) (4,~4) (4,~6) (4,~8) (4,~10) (4,~12) (4,~14) (4, 16) (4, 18) (4, 20) (5, 3) (5, 7) (5,~11) (5,~15) (5,~19) (5,~21) (5,~23) (6, 0) (6, 2) (6, 4) (6, 6) (6, 8) (6,~10) (6,~12) (6,~14) (6,~16) (6,~18) (6,~22) (7,~1) (7, 5) (7, 9) (7, 13) (7, 17) (7,~19) (7,~21) (7,~23) (8,~0) (8,~2) (8,~4) (8,~6) (8, 8) (8, 10) (8, 12) (8, 14) (8,~16) (8,~20) (9,~3) (9,~7) (9,~11) (9,~15) (9,~17) (9, 19) (9, 21) (9, 23) (10, 0) (10,~2) (10,~4) (10,~6) (10,~8) (10,~10) (10,~12) (10,~14) (10,~18) (10, 22) (11,~1) (11,~5) (11,~9) (11,~13) (11,~15) (11,~17) (11,~19) (11,~21) (11, 23) (12,~0) (12,~2) (12,~4) (12,~6) (12,~8) (12,~10) (12,~12) (12,~16) (12,~20) (13,~3) (13,~7) (13,~11) (13,~13) (13,~15) (13,~17) (13,~19) (13,~21) (13,~23) (14,~0) (14,~2) (14,~4) (14,~6) (14,~8) (14,~10) (14,~14) (14,~18) (14,~22) (15,~1) (15,~5) (15,~9) (15,~11) (15,~13) (15,~15) (15,~17) (15,~19) (15,~21) (15,~23) (16,~0) (16,~2) (16,~4) (16,~6) (16,~8) (16,~12) (16,~16) (16, 20) (17, 3) (17, 7) (17,~9) (17,~11) (17,~13) (17,~15) (17,~17) (17,~19) (17,~21) (17, 23) (18,~0) (18,~2) (18,~4) (18,~6) (18,~10) (18,~14) (18,~18) (18,~22) (19,~1) (19,~5) (19,~7) (19,~9) (19,~11) (19,~13) (19,~15) (19,~17) (19,~19) (19,~21) (19,~23) (20,~0) (20,~2) (20,~4) (20,~8) (20,~12) (20,~16) (20, 20) (21,~3) (21,~5) (21,~7) (21,~9) (21,~11) (21,~13) (21,~15) (21,~17) (21, 19) (21, 21) (21, 23) (22, 0) (22,~2) (22,~6) (22,~10) (22,~14) (22,~18) (22, 22) (23, 1) (23, 3) (23, 5) (23,~7) (23,~9) (23,~11) (23,~13) (23,~15) (23, 17) (23, 19) (23, 21) (23, 23)}

\msk

\opti{25}{242}{0.3872}{(0, 0) (0, 2) (0, 4) (0, 6) (0, 8) (0, 10) (0, 12) (0, 14) (0, 16) (0,~18) (0, 20) (0, 22) (0, 24) (1, 1) (1, 5) (1, 9) (1, 11) (1, 15) (1, 19) (1,~23) (2,~0) (2,~2) (2, 4) (2, 6) (2, 8) (2, 12) (2, 14) (2, 16) (2, 18) (2,~20) (2,~22) (2,~24) (3,~3) (3,~7) (3, 9) (3, 11) (3, 13) (3, 17) (3, 21) (4, 0) (4,~2) (4,~4) (4, 6) (4,~10) (4,~14) (4, 16) (4, 18) (4, 20) (4, 22) (4, 24) (5, 1) (5,~5) (5,~7) (5,~9) (5,~11) (5,~13) (5, 15) (5, 19) (5, 23) (6, 0) (6, 2) (6,~4) (6,~8) (6,~12) (6,~16) (6,~18) (6,~20) (6, 22) (6, 24) (7, 3) (7, 5) (7, 7) (7,~9) (7,~11) (7,~13) (7,~15) (7,~17) (7,~21) (8, 0) (8, 2) (8, 6) (8, 10) (8, 14) (8,~18) (8,~20) (8,~22) (8,~24) (9, 1) (9,~3) (9,~5) (9, 7) (9, 9) (9, 11) (9, 13) (9,~15) (9,~17) (9,~19) (9,~23) (10, 0) (10, 4) (10, 8) (10, 12) (10, 16) (10, 20) (10,~22) (10,~24) (11,~1) (11,~3) (11,~5) (11,~7) (11, 9) (11, 11) (11, 13) (11, 15) (11, 17) (11,~19) (11,~21) (12,~0) (12,~2) (12, 6) (12, 10) (12, 14) (12, 18) (12, 22) (12,~24) (13,~3) (13,~5) (13,~7) (13,~9) (13, 11) (13, 13) (13, 15) (13, 17) (13, 19) (13,~21) (13,~23) (14, 0) (14, 2) (14,~4) (14, 8) (14, 12) (14, 16) (14, 20) (14,~24) (15,~1) (15,~5) (15,~7) (15,~9) (15, 11) (15, 13) (15, 15) (15, 17) (15,~19) (15,~21) (15,~23) (16,~0) (16,~2) (16,~4) (16, 6) (16, 10) (16, 14) (16,~18) (16,~22) (16,~24) (17,~3) (17,~7) (17,~9) (17, 11) (17, 13) (17, 15) (17,~17) (17,~19) (17,~21) (18,~0) (18,~2) (18,~4) (18,~6) (18, 8) (18, 12) (18,~16) (18,~20) (18,~22) (18,~24) (19,~1) (19,~5) (19,~9) (19, 11) (19, 13) (19, 15) (19, 17) (19, 19) (19,~23) (20,~0) (20,~2) (20,~4) (20,~6) (20, 8) (20, 10) (20, 14) (20, 18) (20,~20) (20,~22) (20,~24) (21,~3) (21,~7) (21, 11) (21, 13) (21, 15) (21, 17) (21,~21) (22, 0) (22,~2) (22,~4) (22,~6) (22, 8) (22, 10) (22, 12) (22, 16) (22, 18) (22,~20) (22,~22) (22,~24) (23, 1) (23, 5) (23, 9) (23, 13) (23, 15) (23, 19) (23,~23) (24,~0) (24,~2) (24,~4) (24, 6) (24, 8) (24, 10) (24, 12) (24, 14) (24, 16) (24,~18) (24,~20) (24,~22) (24, 24)}

\msk

\opti{26}{260}{0.38461538461538464}{(0, 0) (0, 2) (0, 4) (0, 6) (0, 8) (0, 10) (0, 12) (0,~14) (0,~16) (0,~18) (0,~20) (0,~22) (0, 24) (1, 1) (1, 5) (1, 9) (1, 13) (1, 17) (1,~21) (1,~25) (2,~0) (2,~2) (2,~4) (2, 6) (2, 8) (2, 10) (2, 12) (2, 14) (2, 16) (2, 18) (2, 20) (2, 22) (2, 24) (3, 3) (3, 7) (3, 11) (3, 15) (3, 19) (3, 23) (3, 25) (4, 0) (4,~2) (4,~4) (4,~6) (4,~8) (4,~10) (4, 12) (4, 14) (4, 16) (4, 18) (4, 20) (4, 22) (5, 1) (5, 5) (5, 9) (5, 13) (5, 17) (5, 21) (5, 23) (5, 25) (6, 0) (6, 2) (6, 4) (6,~6) (6,~8) (6,~10) (6,~12) (6,~14) (6, 16) (6, 18) (6, 20) (6, 24) (7, 3) (7,~7) (7,~11) (7,~15) (7,~19) (7,~21) (7, 23) (7, 25) (8, 0) (8, 2) (8, 4) (8, 6) (8, 8) (8, 10) (8, 12) (8, 14) (8, 16) (8, 18) (8, 22) (9, 1) (9, 5) (9, 9) (9, 13) (9,~17) (9,~19) (9,~21) (9,~23) (9,~25) (10, 0) (10, 2) (10, 4) (10, 6) (10, 8) (10,~10) (10,~12) (10,~14) (10,~16) (10,~20) (10, 24) (11, 3) (11, 7) (11, 11) (11, 15) (11,~17) (11,~19) (11,~21) (11,~23) (11,~25) (12, 0) (12, 2) (12, 4) (12, 6) (12,~8) (12,~10) (12,~12) (12,~14) (12,~18) (12, 22) (13, 1) (13, 5) (13,~9) (13,~13) (13,~15) (13,~17) (13,~19) (13,~21) (13,~23) (13, 25) (14, 0) (14,~2) (14,~4) (14,~6) (14,~8) (14,~10) (14,~12) (14,~16) (14,~20) (14,~24) (15,~3) (15,~7) (15,~11) (15,~13) (15,~15) (15,~17) (15,~19) (15,~21) (15,~23) (15,~25) (16, 0) (16, 2) (16,~4) (16,~6) (16,~8) (16,~10) (16,~14) (16,~18) (16,~22) (17,~1) (17,~5) (17,~9) (17,~11) (17,~13) (17,~15) (17,~17) (17,~19) (17,~21) (17,~23) (17, 25) (18, 0) (18, 2) (18, 4) (18,~6) (18,~8) (18,~12) (18,~16) (18,~20) (18,~24) (19,~3) (19,~7) (19,~9) (19, 11) (19, 13) (19, 15) (19, 17) (19, 19) (19,~21) (19,~23) (19,~25) (20,~0) (20,~2) (20,~4) (20,~6) (20,~10) (20,~14) (20,~18) (20,~22) (21,~1) (2,~ 5) (21,~7) (21, 9) (21, 11) (21, 13) (21, 15) (21, 17) (21,~19) (21,~21) (21,~23) (21,~25) (22,~0) (22, 2) (22, 4) (22, 8) (22,~12) (22,~16) (22,~20) (22,~24) (23,~3) (23, 5) (23, 7) (23, 9) (23, 11) (23, 13) (23,~15) (23,~17) (23,~19) (23,~21) (23,~23) (23,~25) (24,~0) (24, 2) (24, 6) (24,~10) (24,~14) (24,~18) (24,~22) (25,~1) (25,~3) (25,~5) (25, 7) (25, 9) (25,~11) (25,~13) (25,~15) (25,~17) (25,~19) (25,~21) (25, 23) (25, 25)}

\end{document}